\documentclass[final,onefignum,onetabnum,letterpaper]{siamart250211}

\usepackage{soul}
\usepackage{lipsum}
\usepackage{amsfonts}
\usepackage{epstopdf}
\usepackage{epsfig}
\ifpdf
  \DeclareGraphicsExtensions{.eps,.pdf,.png,.jpg}
\else
  \DeclareGraphicsExtensions{.eps}
\fi

\usepackage{datetime}
\newdateformat{monthyeardate}{%
	\monthname[\THEMONTH] \THEDAY, \THEYEAR}

\usepackage{academicons}
\usepackage{xcolor}

\usepackage{amsmath}
\usepackage{amssymb}
\usepackage{commath}
\usepackage{mathtools}
\usepackage{bbm}
\usepackage{bm}
\usepackage{upgreek}
\usepackage{float}
\allowdisplaybreaks
\usepackage{import}
\usepackage{color}
\usepackage{graphicx}
\usepackage[small]{caption}
\usepackage{subcaption}
\usepackage{tabularray}
\UseTblrLibrary{amsmath}

\usepackage{relsize}
\usepackage{adjustbox}
\usepackage{algorithm}
\usepackage{algpseudocode}
\usepackage{booktabs}
\usepackage{tikz}
\usepackage{pifont}
\usetikzlibrary{bayesnet} 

\usepackage{enumitem}
\setlist[enumerate]{leftmargin=.5in}
\setlist[itemize]{leftmargin=.5in}

\newsiamremark{remark}{Remark}
\crefname{hypothesis}{Hypothesis}{Hypotheses}

\headers{Entropy-stable learned Hyperbolic Conservation Laws}{Lizuo Liu, Lu Zhang and Anne Gelb}

\title{Parametric Hyperbolic Conservation Laws: A Unified Framework for Conservation, Entropy Stability, and Hyperbolicity
\thanks{
\monthyeardate\today 
}}

\author{
Lizuo Liu\thanks{Department of Mathematics, Dartmouth College, Hanover, NH 03755, USA 
(\email{Lizuo.Liu@dartmouth.edu}, 
\email{Anne.E.Gelb@dartmouth.edu}
)}
\and
Lu Zhang\thanks{Department of Computational Applied Mathematics and Operations Research, and Ken Kennedy Institute, Rice University, Houston, TX 77005, USA}
(\email{lz82@rice.edu})
\and 
Anne Gelb\footnotemark[2]
}

\usepackage{amsopn}

\usepackage{comment}

\usepackage[normalem]{ulem}

\newcommand{\bmu}{{\boldsymbol\mu}}
\newcommand{\btheta}{{\boldsymbol\theta}}

\newcommand{\bu}{\mathbf{u}}

\usepackage[normalem]{ulem}

\newtheorem{defi}{Definition}[section]

\numberwithin{equation}{section}
\numberwithin{figure}{section}
\numberwithin{table}{section}

\allowdisplaybreaks

\patchcmd\newpage{\vfil}{}{}{}
\begin{document}

\maketitle

\begin{abstract}
	We propose a parametric hyperbolic conservation law (SymCLaw) for learning hyperbolic systems directly from data while ensuring conservation, entropy stability, and hyperbolicity by design. Unlike existing approaches that typically enforce only conservation or rely on prior knowledge of the governing equations, our method parameterizes the flux functions in a form that guarantees real eigenvalues and complete eigenvectors of the flux Jacobian, thereby preserving hyperbolicity. At the same time, we embed entropy-stable design principles by jointly learning a convex entropy function and its associated flux potential, ensuring entropy dissipation and the selection of physically admissible weak solutions. A corresponding entropy-stable numerical flux scheme provides compatibility with standard discretizations, allowing seamless integration into classical solvers. Numerical experiments on benchmark problems, including Burgers’, shallow water, Euler, and KPP equations, demonstrate that SymCLaw generalizes to unseen initial conditions, maintains stability under noisy training data, and achieves accurate long-time predictions, highlighting its potential as a principled foundation for data-driven modeling of hyperbolic conservation laws.
\end{abstract}

\begin{keywords}
   	Parametric conservation laws, conservation, hyperbolicity, entropy stability, data-driven modeling
\end{keywords}

\begin{AMS}
        65M08, 68T07, 65M22, 65M32, 65D25
\end{AMS}



\section{Introduction}\label{sec: introduction}
Conservation laws provide the mathematical foundation for modeling a wide range of physical systems, including fluid dynamics \cite{Euler,majda2012compressible}, traffic flow \cite{delle2014scalar}, magnetohydrodynamics \cite{webb2018magnetohydrodynamics}, geophysics \cite{shepherd1990symmetries}, and many more. By encoding fundamental invariants such as mass, momentum, and energy,  equations governed by conservation laws naturally capture the nonlinear wave phenomena observed in practice. Their broad utility in describing large-scale flow patterns and wave phenomena also makes them fundamental to modeling weather and climate systems \cite{Euler}. A particularly relevant example involves iceberg and sea ice dynamics, which affect both oceanic circulation and atmospheric processes; here, hyperbolic PDEs are used to model iceberg motion \cite{iceberg} and the elastic or plastic deformation of sea ice \cite{iceberg, Girard11, Elastic-PlasticSeaIce}. Yet in many emerging scientific and engineering applications,  including sea ice dynamics modeling \cite{iceberg, Girard11, Elastic-PlasticSeaIce}, the explicit form of the underlying conservation laws is either unknown or only partially specified through empirical models. In such cases, abundant trajectory data from high-fidelity simulations or experiments often serve as the primary source for uncovering or approximating the underlying conservation laws. This need has motivated a growing body of work on data-driven approaches to learning conservation laws. 


A central challenge in learning conservation laws is ensuring that the inferred dynamics are not only consistent with observed data but also physically admissible in a way that preserves the essential structure of the equations. Standard data-driven approaches, such as black-box regression models, may achieve good fits to observed trajectories, but without careful design they risk generating solutions that exhibit spurious oscillations, unphysical behaviors, or ill-posed dynamics outside the training regime. Several recent approaches have addressed these issues  by emphasizing the enforcement of  fundamental conservation  properties within the learned models. For example, Godunov--Riemann informed neural networks \cite{Dimitrios2025gorinn} embed Riemann-solver structures to learn physical flux functions, while the framework in \cite{Seung2025tvdnet} incorporates neural network closures inspired by total variation diminishing methods to recover dynamics. Other notable contributions include the conservative flux form neural network (CFN)~\cite{chencfn}, symbolic frameworks that mimic finite volume schemes \cite{li2023flux, li2024alternating}, and RoeNet \cite{RoeNet}, which integrates the classical Roe solver into a neural architecture. Together, these efforts have laid important groundwork for data-driven modeling of conservation laws by explicitly embedding conservation principles into the learning process. 

However, for stable and accurate simulation of nonlinear conservation laws, conservation alone is not sufficient: the emergence of shocks necessitates additional admissibility criteria, such as entropy stability, to ensure physically meaningful weak solutions. This observation has motivated a subsequent line of research focused on extending data-driven approaches beyond conservation to incorporate entropy-based stability principles. For example, \cite{liu2024entropy} introduced the entropy-stable CFN (ESCFN), which augments the CFN framework with entropy stability through slope-limiting mechanisms embedded in the neural network architecture. Building on this, \cite{liu2025neural} further proposed to learn not only the fluxes but also the entropy function itself, thereby enabling the model to select physically admissible solutions among the family of weak solutions. These advances mark an important step toward embedding deeper structural properties of conservation laws into machine-learned models. Nonetheless, a fundamental challenge remains unresolved: ensuring that the learned system is \emph{hyperbolic}.

Hyperbolicity, the requirement that the flux Jacobian has real eigenvalues and a complete set of eigenvectors, is a property essential for well-posedness and wave-propagation dynamics. It ensures well-posedness of the Cauchy problem, stability of numerical approximations, and the correct finite-speed propagation of information. Without hyperbolicity, the dynamics can lose predictive meaning, as wave propagation may become ill-defined or unstable. At the same time, hyperbolicity is closely connected to entropy: the existence of a strictly convex entropy-entropy flux pair, which provides the admissibility criterion for shocks, also serves to symmetrize the system and is guaranteed only in the hyperbolic setting (see \Cref{thm: equivalence}). This structural interplay highlights that hyperbolicity is not merely a technical condition but a critical requirement for both the mathematical soundness and physical fidelity of conservation laws. For comprehensive discussions on this topic, see \cite{ Barth1999NumericalMethodsGasDynamic, dafermos2005hyberbolic, Tadmor1987EntropyStable, tadmor2016entropy}.

In this work, we introduce a parametric hyperbolic conservation law that, to our knowledge, is the first to simultaneously enforce \emph{conservation}, \emph{entropy stability}, and \emph{hyperbolicity} within a unified learning paradigm. The central idea is to parameterize the fluxes in a form that guarantees the flux Jacobian matrices admit real eigenvalues and a complete set of eigenvectors, thereby ensuring hyperbolicity. At the same time, we learn a strictly convex entropy function, modeled via an input-convex neural network (ICNN), together with its associated flux potential. This design leverages the fact that a strictly convex entropy symmetrizes the system of conservation laws (see \Cref{thm: equivalence}), so that hyperbolicity is preserved by construction while entropy stability is naturally embedded. Building on this structurally guaranteed formulation, we further propose a differentiable entropy-stable numerical scheme that integrates seamlessly with classical solvers, providing a principled bridge between data-driven modeling and established numerical methods. This dual guarantee, consistency with data and adherence to fundamental mathematical structure, not only preserves the essential properties of conservation laws but also enhances robustness when extrapolating beyond the training regime. Through a series of benchmark problems, including Burgers’, shallow water, Euler, and KPP\footnote{The KPP equations are named for  Kurganov, Petrova, and Popov in \cite{kpp}.} equations, we demonstrate that our framework achieves accurate long-time predictions and maintains stability even under noisy training data, outperforming existing approaches that enforce only conservation or entropy stability.

The rest of this paper is organized as follows. In \Cref{sec: preliminaries} we review the key theoretical and numerical foundations relevant to hyperbolic systems of conservation laws. \Cref{sec: parametric conservation laws} introduces parametric hyperbolic conservation laws (SymCLaw) and the corresponding entropy-stable flux schemes for learning hyperbolic conservation laws directly from solution trajectories. The experimental setup and evaluation metrics are described in \Cref{sec: training design}. In \Cref{sec: numerics}, we present a series of numerical examples that demonstrate the effectiveness and robustness of the proposed SymCLaw  method. Finally, we conclude with a summary and briefly discuss directions for future research in \Cref{sec: conclusion}.

\section{Preliminaries}
\label{sec: preliminaries}
 We begin by reviewing the essential theoretical foundations of hyperbolic systems of conservation laws, which serve as the basis for SymCLaw introduced in \Cref{sec: parametric conservation laws}.

\subsection{Hyperbolic conservation laws}
\label{sec:conservationlaws}
 Our primary objective is to learn hyperbolic conservation laws directly from observed data and to employ the learned models for predicting future dynamics. This approach is particularly important in scenarios where the governing equations are not available in closed form. To ensure that the learned model generates physically relevant solutions, it is essential to identify the key structural properties satisfied by such systems. To this end, we provide a brief review of these properties.
 
 Consider the $d$-dimensional system of conservation laws 
\begin{equation} \label{eq: conservation_law} \frac{\partial \mathbf{u}}{\partial t}  + \sum_{\mathfrak{i}=1}^d \frac{\partial \mathbf{f}_\mathfrak{i}(\mathbf{u})}{\partial x_\mathfrak{i}}  = \mathbf{0}, \quad \mathbf{x} \in \Omega \subset \mathbb{R}^d, \quad t \in (0, T), 
\end{equation} 
subject to appropriate initial and boundary conditions. Here, $\mathbf{u} = [u^1, \ldots, u^p]^\top$ represents the vector of conserved state variables, taking values in a convex set $\mathcal{D} \subset \mathbb{R}^p$. The flux functions $\mathbf{f}_\mathfrak{i} : \mathcal{D} \rightarrow \mathbb{R}^p$ are assumed to be sufficiently smooth. For each spatial direction $1 \leq \mathfrak{i} \leq d$, we define the flux Jacobian matrix as
\begin{equation}\label{eq: flux_jacobi}
A_\mathfrak{i}({\bf u}) := \nabla_{\bf u}\mathbf{f}_\mathfrak{i}(\mathbf{u})=\left\{\frac{\partial f_\mathfrak{i}(\mathbf{u})}{\partial u^\mathfrak{k}}\right\}_{1 \leq \mathfrak{i}, \mathfrak{k} \leq p}.
\end{equation}
The system \eqref{eq: conservation_law} is called \emph{hyperbolic} if for any unit vector $\mathbf{n} = (n_1, n_2, \cdots, n_d)^\top \in \mathbb{R}^d$ and all $\mathbf{u} \in \mathcal{D}$, the matrix $\sum_{\mathfrak{i}=1}^d n_\mathfrak{i} A_\mathfrak{i}(\mathbf{u})$ has $p$ real eigenvalues and a complete set of linearly independent eigenvectors.

\subsubsection{Entropy solution} \label{sec:entropysolution} It is well known that even when the initial data are smooth, solutions to \eqref{eq: conservation_law} may develop discontinuities in finite time. Consequently, solutions of \eqref{eq: conservation_law} must be interpreted in the weak (distributional) sense. However, weak solutions are in general not unique. To enforce uniqueness and ensure physical admissibility, additional constraints are imposed in the form of entropy conditions.
\begin{defi}
Suppose that the domain $\mathcal{D}$ is convex. A convex function $\eta(\mathbf{u}): \mathcal{D} \rightarrow \mathbb{R}$ is called a convex \emph{entropy function} for the system \eqref{eq: conservation_law} if there exist associated entropy fluxes $G_\mathfrak{i}(\mathbf{u}): \mathcal{D}\rightarrow \mathbb{R}$, $1 \leq \mathfrak{i} \leq d$, such that
\[\nabla_{\bf u}\eta(\mathbf{u}) \nabla_{\bf u}\mathbf{f}_\mathfrak{i}(\mathbf{u})=\nabla_{\bf u}G_\mathfrak{i}(\mathbf{u})\quad \forall \mathfrak{i},\]
where the gradient, $\nabla_{\bf u}$, is understood as row vector.
\end{defi}

In regions where the solution is smooth, multiplying \eqref{eq: conservation_law} on the left by $\nabla_{\bf u}\eta(\mathbf{u})$ yields an additional conservation law
\[\frac{\partial \eta(\mathbf{u})}{\partial t}+\sum_{\mathfrak{i}=1}^d\frac{\partial G_\mathfrak{i}(\mathbf{u})}{\partial x_\mathfrak{i}}=0.\]
Across discontinuities, however, entropy is required to dissipate, which leads to the entropy inequality condition
\begin{equation} \label{eq:entropy_cond} 
\frac{\partial \eta(\mathbf{u})}{\partial t} + \sum_{\mathfrak{i}=1}^d \frac{\partial G_\mathfrak{i}(\mathbf{u})}{\partial x_\mathfrak{i}} \leq 0,
\end{equation}
and selects the physically admissible weak solution, commonly referred to as the \emph{entropy solution}.

\subsubsection{Symmetrization} \label{sec:symmetrization}
When the entropy function $\eta({\bf u})$ is strictly convex, the transformation $\mathbf{v} = (\nabla_{\bf u}\eta(\mathbf{u}))^\top$ defines a one-to-one change of variables from $\mathbf{u}$ to the so-called \emph{entropy variables} $\mathbf{v}$. We write the fluxes in terms of these variables as $\mathbf{g}_\mathfrak{i}(\mathbf{v}) := \mathbf{f}_\mathfrak{i}(\mathbf{u}(\mathbf{v}))$, thereby allowing us to recast the conservation law \eqref{eq: conservation_law} in the symmetric form
\begin{equation}\label{eq: symmetric}
\big(\nabla_{\bf v}\mathbf{u}(\mathbf{v}) \big)\frac{\partial \mathbf{v}}{\partial t}+\sum_{\mathfrak{i}=1}^d \nabla_{\bf v}\mathbf{g}_\mathfrak{i}(\mathbf{v}) \frac{\partial \mathbf{v}}{\partial x_\mathfrak{i}}=0.
\end{equation}
Due to the strict convexity of $\eta({\bf u})$, the matrix $\nabla_{\bf v}\mathbf{u}(\mathbf{v}) = \left( \mathbb{H}_{\bf u}\eta(\mathbf{u}) \right)^{-1}$ is symmetric and positive-definite. The following theorem, adapted from \cite{godlewski2013numerical} (Theorems 5.1 and 5.2), establishes the equivalence between the existence of a strictly convex entropy function and the symmetry structure of the system. 

\begin{theorem}\label{thm: equivalence}
   Let $\eta({\bf u})$ be a strictly convex function. Then $\eta({\bf u})$ serves as an entropy function for \eqref{eq: conservation_law} if and only if the matrix $\nabla_{\bf v}\mathbf{u}(\mathbf{v})$ is symmetric positive-definite and $\nabla_{\bf v}\mathbf{g}_\mathfrak{i}(\mathbf{v})$ is symmetric for each $1 \leq \mathfrak{i} \leq d$. The system \eqref{eq: symmetric} is then called the \emph{symmetrized form} of \eqref{eq: conservation_law}. Furthermore, for any unit vector $\mathbf{n} = (n_1, n_2, \cdots, n_d)^\top \in \mathbb{R}^d$, the matrix
   \[A(\mathbf{u}, \mathbf{n}):=\sum_{\mathfrak{i}=1}^d n_\mathfrak{i} \nabla_{\bf u}\mathbf{f}_\mathfrak{i}(\mathbf{u})=\sum_{\mathfrak{i}=1}^d n_\frak{i} \nabla_{\bf v}\mathbf{g}_\mathfrak{i}(\mathbf{v}) \nabla_{\bf u}\mathbf{v}(\mathbf{u})\] is similar  to the symmetric matrix
   \[\big(\nabla_{\bf u}\mathbf{v}(\mathbf{u})\big)^{\frac{1}{2}}\Big(\sum_{\mathfrak{i}=1}^d n_\mathfrak{i} \nabla_{\bf v}\mathbf{g}_\mathfrak{i}(\mathbf{v})\Big) \big(\nabla_{\bf u}\mathbf{v}(\mathbf{u})\big)^{\frac{1}{2}}.\] 
   Consequently, the existence of a strictly convex entropy function implies that the system \eqref{eq: conservation_law} is hyperbolic.
\end{theorem}

This symmetrization result plays a fundamental role in both the theoretical analysis and computational treatment of classical hyperbolic conservation laws. Building on this foundation, we introduce what we will call the {\em SymCLaw} method to predict  the dynamics of an {\em unknown} system of hyperbolic conservation laws that is fully derived from available data trajectories over a limited time domain.  Importantly, the SymClaw method is able to  preserve essential physical properties and does not require oracle knowledge of the solution space.

\section{Parametric hyperbolic conservation laws (SymCLaw)}\label{sec: parametric conservation laws}
 In many real-world applications, the equations governing dynamics are inaccessible and can only possibly be inferred from available data trajectories.  For example, one might consider how measurements from fluid dynamics, traffic flow, or wave propagation might infer their corresponding PDE models. In this investigation we consider that  these phenomena are governed by conservation principles and then ask the fundamental question: \emph{Can we reconstruct the underlying physical laws directly from limited data trajectories, assuming they adhere to a hyperbolic conservation structure?} Motivated by \Cref{thm: equivalence}, we emphasize that the presence of a strictly convex entropy function, together with the symmetry of the flux Jacobian in entropy variables, guarantees the hyperbolicity of the corresponding conservation law, thereby providing a principled foundation for such data-driven reconstruction.

Inspired by the methodology developed in \cite{liu2025neural}, here we propose a parametric formulation guided by a learnable strictly convex entropy function, denoted by \(\eta_{\btheta}({\bf u})\), where $\btheta$ represents the parameters of an input convex neural network. For clarity of presentation, we defer the architectural details of the neural networks to  \Cref{sec:networkdetails}. According to \Cref{thm: equivalence}, $\eta_{\btheta}({\bf u})$ qualifies as an entropy function associated with a conservation law if and only if the flux Jacobian, expressed in terms of the parametrized entropy variable, 
\[{\bf v} = (\nabla_{\bf u}\eta_\btheta({\bf u}))^{\top},\] 
is symmetric. This condition holds because $\nabla_{\bf v}({\bf u}({\bf v})) = (\mathbb{H}_{{\bf u}}\eta_\btheta({\bf u}))^{-1}$ is symmetric by the construction of the convex entropy function. To enforce the hyperbolicity, we then model the flux functions ${\bf g}_{\mathfrak{i}}({\bf v})$ in \eqref{eq: symmetric} using gradient neural networks: 
\begin{equation}\label{eq:gradient_NN}
{\bf g}_{\mathfrak{i}}({\bf v}) = \nabla_{\bf v}\phi_{\bmu, \mathfrak{i}}({\bf v}), \quad \mathfrak{i} = 1,2,\cdots, d,
\end{equation}
where \(\phi_{\bmu,\mathfrak{i}}({\bf v})\) is a scalar-output neural network with parameters $\bmu$. By construction, the symmetry of ${\bf g}_{\mathfrak{i}}({\bf v})$ is automatically guaranteed, since the Hessian of any scalar function is symmetric. In terms of hyperbolic conservation laws, $\phi_{\bmu, \mathfrak{i}}({\bf v})$ is called the entropy flux potential since its gradient produces the flux function ${\bf g}_{\mathfrak{i}}({\bf v})$.

With these ingredients in place, we are now ready to introduce what we will refer to as the SymClaw method:
\begin{equation}\label{eq: paramatric HCLaw}
 \frac{\partial \mathbf{u}}{\partial t}+\sum_{\mathfrak{i}=1}^d  \frac{\partial }{\partial x_{\mathfrak{i}}} \underbrace{\nabla_{\bf v}\phi_{\bmu, \mathfrak{i}}\big(\overbrace{\nabla_{\bf u}\eta_{\btheta}({\bf u})}^{\bf v^\top}\big)}_{{\bf f}_{\mathfrak{i}}^{\btheta,\bmu}({\bf u})}={\bf 0}.
\end{equation}
In particular, the flux Jacobian in terms of the state variable $\bf u$ takes the form
\[
\nabla_{\bf u} {\bf f}^{\btheta,\bmu}_\mathfrak{i}({\bf u}) = \left(\mathbb{H}_{\bf v} \phi_{\bmu, \mathfrak{i}}({\bf v})\right)\left(\mathbb{H}_{\bf u}\eta_{\btheta}({\bf u})\right),
\]
which are the product of two Hessians. By construction of network structures we have that \(\mathcal{A}_\mathfrak{i} := \mathbb{H}_{\bf v} \phi_{\bmu, \mathfrak{i}}({\bf v})\) is symmetric, while \( \mathcal{B} := \mathbb{H}_{\bf u}\eta_{\btheta}({\bf u})\) is symmetric positive definite. Consequently, by \Cref{thm: equivalence}, the matrix $\mathcal{A}_\mathfrak{i}\mathcal{B}$ is similar to the symmetric matrix $\mathcal{B}^{\frac{1}{2}}\mathcal{A}_\mathfrak{i}\mathcal{B}^{\frac{1}{2}}$, and therefore has only real eigenvalues. It follows that  \eqref{eq: paramatric HCLaw} forms a system of (parametric) hyperbolic conservation laws, with $\eta_\btheta({\bf u})$ serving as an entropy function, and $\phi_{\bmu,\mathfrak{i}}({\bf v}), \mathfrak{i} = 1,\cdots,d$, serving as a corresponding entropy flux potential.

We are now left with the fundamental question: \emph{How should the parameter set $(\btheta,\bmu)$  in \cref{eq: paramatric HCLaw} be determined?} While there may be several ways to approach this, here we propose to identify these parameters in the discrete setting, that is, by first discretizing  \eqref{eq: paramatric HCLaw} to obtain its discrete counterpart, and then to use the available data trajectories to calculate $(\btheta,\bmu)$.  The rationale for working in the discrete setting rather than directly in the continuous formulation is twofold. First, observed data, whether from experiments or high-fidelity simulations, are inherently discrete in time and space, making direct comparison to the continuous model both ill-posed and numerically unstable. Second, identifying parameters at the continuous PDE level requires handling functional equations involving fluxes and entropy pairs, which can be analytically intractable and sensitive to noise. By contrast, the discrete SymCLaw naturally aligns with the format of the data, provides a stable optimization framework, and allows us to leverage well-established numerical schemes. It is worth noting that the proposed approach is not tied to any specific numerical scheme, and  in principle, one could employ finite volume, finite element, or discontinuous Galerkin methods, among others. Here we adopt a finite volume discretization due to its simplicity. 

\subsection{Finite volume methods for conservation laws}\label{sec:FiniteVolume}  We  begin by briefly reviewing  classical finite volume schemes for hyperbolic conservation laws, with particular emphasis on the construction of entropy-stable numerical fluxes. Such fluxes play a crucial role in ensuring robustness, especially in problems where discontinuities or shocks dominate the solution behavior, making entropy stability an indispensable property of reliable numerical schemes.  This background will provide the means for determining the parameter set $(\btheta,\bmu)$ in \cref{eq: paramatric HCLaw} from the corresponding observed data trajectories.

For illustrative purposes, we consider \eqref{eq: conservation_law} with $d=1$,\footnote{We emphasize that although for ease of presentation we use the 1D set up, the proposed framework naturally extends to multidimensional systems.} given by
\begin{equation}\label{eq:conservation_law_1d}
\frac{\partial \mathbf{u}}{\partial t} + \frac{\partial \mathbf{f}(\mathbf{u})}{\partial x} = 0, \quad x \in \Omega = (a,b), \quad t \in (0,T).
\end{equation} 
Here, ${\bf u}\in \mathbb{R}^p$ represents the conserved states, and ${\bf f}({\bf u}) \in \mathbb{R}^p$ is the corresponding flux function. For the SymCLaw equation \cref{eq: paramatric HCLaw}, ${\bf f}({\bf u})$ is replaced by ${\bf f}^{\btheta,\bmu}({\bf u})$. 

To numerically solve \eqref{eq:conservation_law_1d}, we partition the computational domain $\Omega$ into uniform mesh with grid points $\{x_j\}_{j=0}^n$, where $x_j= j\Delta x$ and $\Delta x = \frac{b-a}{n}.$ The solution is approximated by cell averages ${\bf u}_j(t) = [u_j^1(t), \cdots, u_j^p(t)]^\top$ over the cell $I_j := (x_j - \frac{\Delta x}{2}, x_j + \frac{\Delta x}{2})$ with 
\begin{equation}\label{eq: cell_average}
u_j^{\mathfrak{i}}(t)=\int_{I_j} u^{\mathfrak{i}}(x, t) d x, \quad j=1, \cdots, n-1; \ \mathfrak{i} = 1,\cdots,p,
\end{equation}
and appropriate boundary conditions are enforced at the endpoints, denoted by $u_0^{\mathfrak{i}}$ and $u_n^{\mathfrak{i}}$. For notational simplicity, we use ${\bf u} = [u^1, \cdots, u^p]^\top$ to denote the continuous solution and ${\bf u}_j = [u_j^1, \cdots, u_j^p]^\top$ to denote its cell average over $I_j$. This convention will be adopted throughout the paper unless otherwise specified, and any potential ambiguity will be explicitly clarified.

Within this framework, let $v_j^{\mathfrak{i}}(t), 1\leq {\mathfrak{i}}\leq p$, denote the approximation to the entropy variable $v^{\mathfrak{i}}(x,t)$ averaged over $I_j$ using \eqref{eq: cell_average}, and define ${\bf v}_j = [v_j^1, \cdots, v_j^p]^\top$. We then have the following result, adapted from \cite{Tadmor1987EntropyStable, tadmor2016entropy}, where three-point entropy-stable finite-volume methods were originally established.

\begin{theorem}\label{thm:entropy-stable}
Let $(\eta({\bf u}), G({\bf u}))$ be an entropy pair associated with the conservation law~\eqref{eq:conservation_law_1d}, and let ${\bf v} = \nabla_{\bf u}\eta({\bf u})$ denote the corresponding entropy variable.  
Define ${\bf g}({\bf v}) := {\bf f}({\bf u}({\bf v}))$ and the entropy flux potential
\begin{equation}\label{eq: entropy_potential}
\phi({\bf v}) := {\bf v}^\top {\bf g}({\bf v}) - G({\bf u}({\bf v})).
\end{equation}
Let ${\bf f}^\ast_{j+1/2}$ be any entropy-conservative numerical flux satisfying
\begin{equation}\label{eq:entropy_conservative condition}
[[{\bf v}]]_{j+1/2}^\top {\bf f}^\ast_{j+1/2} = [[\phi({\bf v})]]_{j+1/2},
\end{equation}
where $[[{\bf v}]]_{j+1/2} := {\bf v}_{j+1/2}^+  - {\bf v}_{j+1/2}^-$ and $[[\phi({\bf v})]]_{j+1/2} := \phi({\bf v}_{j+1/2}^+) - \phi({\bf v}_{j-1/2}^-)$, with ${\bf v}_{j+1/2}^-$ and ${\bf v}_{j+1/2}^+$ denoting the reconstructed entropy variables at the left and right sides of the cell interface $x_{j+1/2}$, respectively. Then, the semi-discrete finite-volume scheme
\[
\frac{d}{dt}\,{\bf u}_j(t)
+
\frac{1}{\Delta x}
\left(
{\bf f}^\ast_{j+1/2} - {\bf f}^\ast_{j-1/2}
\right)
= {\bf 0},
\]
is \emph{entropy-conservative} in the sense that
\[
\frac{d}{dt}\,\eta({\bf u}_j(t))
+
\frac{1}{\Delta x}
\left(
G_{j+1/2} - G_{j-1/2}
\right)
= 0,
\]
where $G_{j+1/2}$ is a numerical entropy flux consistent with the entropy flux potential $\phi({\bf v})$ defined in \eqref{eq: entropy_potential}. Furthermore, define the \emph{entropy-stable numerical flux}
\begin{equation}\label{eq: entropy stable flux}
\widehat{\bf f}_{j+1/2}
=
{\bf f}^\ast_{j+1/2}
-
\tfrac{1}{2}\,D_{j+1/2}\,[[{\bf v}]]_{j+1/2},
\end{equation}
where $D_{j+1/2} \succeq 0$ is a symmetric positive semidefinite matrix. Then, the semi-discrete finite-volume scheme
\[
\frac{d}{dt}\,{\bf u}_j(t)
=
-\frac{1}{\Delta x}
\left(
\widehat{\bf f}_{j+1/2} - \widehat{\bf f}_{j-1/2}
\right)
\]
is \emph{entropy-stable}, in the sense that
\[
\frac{d}{dt}\,\eta({\bf u}_j(t))
+
\frac{1}{\Delta x}
\left(
\widehat G_{j+1/2} - \widehat G_{j-1/2}
\right)
\le 0,
\]
with numerical entropy flux
\[
\widehat G_{j+1/2}
=
G_{j+1/2}
+
\tfrac{1}{2}\,\overline{\bf v}_{j+1/2}^\top
D_{j+1/2}\,[[{\bf v}]]_{j+1/2},
\quad \text{where} \quad
\overline{\bf v}_{j+1/2} := \tfrac{1}{2}\,({\bf v}_{j+1/2}^- + {\bf v}_{j+1/2}^+).
\]
\end{theorem}

We note that the construction of entropy-stable schemes follows the classical  framework in \cite{tadmor2016entropy}, in which an entropy-conservative flux ${\bf f}_{j+1/2}^\ast$ is augmented with a dissipation term through a symmetric positive semidefinite matrix $D_{j+1/2}$. Various choices of entropy-conservative flux ${\bf f}_{j+1/2}^\ast$ are possible, depending on the desired accuracy and the structure of the conservation law. Here  we construct the entropy conservative flux by projecting the arithmetic average of the physical fluxes onto the direction of the entropy-variable jump, following the discrete entropy-conservation condition \eqref{eq:entropy_conservative condition}, that is,   
\begin{equation}\label{eq: entropy conservative flux}
\begin{aligned}
{\bf f}^\ast_{j+1/2} &= \tfrac{1}{2}\big({\bf f}({\bf u}_{j+1/2}^+)+{\bf f}( {\bf u}_{j+1/2}^-)\big) \\
&+ \frac{[[\phi({\bf v})]]_{j+1/2} - \frac{1}{2}[[{\bf v}]]^\top\big({\bf f}({\bf u}_{j+1/2}^+) + {\bf f}({\bf u}_{j+1/2}^-)\big)}{\|[[{\bf v}]]_{j+1/2}\|^2_2} [[{\bf v}]]_{j+1/2},
\end{aligned}
\end{equation}
which satisfies the discrete entropy-conservation \eqref{eq:entropy_conservative condition} for any pair of states ${\bf u}_{j+1/2}^{\pm}$. The overall order of accuracy of the finite-volume scheme depends on the accuracy of these reconstructed interface states.  
If ${\bf u}_{j+1/2}^\pm$ are obtained as piecewise constant values (i.e., ${\bf u}_j$ and ${\bf u}_{j+1}$), the scheme is formally second-order accurate. However, when ${\bf u}_{j+1/2}^\pm$ is reconstructed using a $p$th-order accurate procedure, the resulting finite-volume method attains $p$th-order accuracy while retaining the exact discrete entropy-conservation property  for smooth  regions. In this work, the left and right interface states ${\bf u}_{j+1/2}^{\pm}$ are reconstructed using the fifth-order Weighted Essentially Non-Oscillatory (WENO5) scheme (\cite{shu2020essentially}). This choice provides high-order accuracy in smooth regions while effectively controlling spurious oscillations near discontinuities, thereby combining entropy stability with nonlinear robustness.


\subsection{Entropy-stable flux for SymCLaw} \label{sec:esfluxsymclaw} 


We now describe in detail how entropy-stable fluxes are encoded within the framework of SymCLaw.  In particular, the  parameters $\btheta$ and $\bmu$ in \eqref{eq: paramatric HCLaw} are determined by minimizing the misfit between the observed data trajectories and those generated by the discrete model, where the model is intentionally solved using an entropy-stable finite volume scheme (see \cref{eq: entropy stable flux}). The derivation of the semidiscretized system of ODEs in \cref{eq: scheme_learned_law} at the end of this section is key to this process.

As stated in \Cref{thm:entropy-stable}, the construction of entropy-stable numerical schemes for hyperbolic conservation laws relies on augmenting entropy-conservative fluxes with appropriate dissipation. This is achieved in practice by adding a numerical viscosity term to an entropy-conservative flux, typically constructed using a symmetric positive semidefinite diffusion matrix together with the jump of entropy variables across adjacent cells (see \eqref{eq: entropy stable flux}). While this general strategy provides a flexible foundation for ensuring entropy stability, practical implementations often employ specific choices that align with the structure of the underlying conservation law. A common approach is to derive the diffusion matrix from the eigendecomposition of the flux Jacobian (see, e.g., \cite{fjordholm2012arbitrarily,FarzadAffordableEntropyConsistentEulerFlux}). Following our recent work \cite{liu2025neural}, we adopt a Rusanov-type entropy-stable flux for \eqref{eq: paramatric HCLaw}:
\begin{equation}
    \widehat{\bf f}^{\btheta,\bmu}_{j+1/2} =  {\bf f}^{\btheta,\bmu,\ast}_{j+1/2} - \tfrac{1}{2} \underbrace{\lambda_{j+1/2}^{\max} \Big(\mathbb{H}_{\bf u}\eta_{\btheta}(\bar{\bf u}_{j+1/2})\Big)^{-1}}_{D_{j+1/2}} [[\underbrace{\nabla_{\bf u}\eta_\btheta({\bf u})}_{{\bf v}^\top}]]_{j+1/2},
\label{eq: Param entropy stable flux}
\end{equation}
where 
\begin{equation}\label{eq:maximum_wave_speed}
\lambda_{j+1/2}^{\max} = \rho(\mathcal{B}^{\frac{1}{2}} \mathcal{A} \mathcal{B}^{\frac{1}{2}}), \quad \text{with} \quad \mathcal{A} = \mathbb{H}_{\bf v} \phi_{\bmu}(\bar{\bf v}_{j+1/2})\  \text{ and } \ \mathcal{B} = \mathbb{H}_{\bf u}\eta_{\btheta}(\bar{\bf u}_{j+1/2}),
\end{equation}
is the local maximum wave speed at cell interface $x_{j+1/2}$. In contrast to \cite{liu2025neural}, where $\lambda_{j+1/2}^{\text{max}}$ is approximated using a neural network, the computation of $\lambda_{j+1/2}^{\text{max}}$ in the present work is both accurate and GPU-efficient, as JAX provides optimized batched implementations of Jacobi iteration methods for small Hermitian matrices. 

We further emphasize that the choice of Rusanov-type entropy-stable fluxes in this study is primarily motivated by considerations of computational efficiency, and that the proposed methodology can be naturally extended to other classes of fluxes. For instance, Roe-type entropy-stable fluxes may be employed, in which the full spectrum of wave speeds and the corresponding eigenvectors are incorporated rather than relying solely on the maximum local wave speed. However, the stable and efficient computation of the complete eigensystem of the flux Jacobian remains challenging and lies beyond the scope of this work, and is thus left for future investigations.

To complete the formulation in \eqref{eq: Param entropy stable flux}, the entropy-conservative flux $\mathbf{f}^{\btheta,\bmu,\ast}_{j+1/2}$ must also be specified. Analogous to  \cref{eq: entropy conservative flux}, here we construct 
\begin{equation}\label{eq:entropy_stable_flux_nn}
\begin{aligned}
{\bf f}_{j+1/2}^{\btheta,\bmu,\ast} &\ = \tfrac{1}{2}\big({\bf f}^{\btheta,\bmu}({\bf u}_{j+1/2}^+)\!+\!{\bf f}^{\btheta,\bmu}( {\bf u}_{j+1/2}^-)\big)  \\
& \ + \frac{[[\phi_{\bmu}(\nabla_{\bf u}\eta_\btheta({\bf u}))]]_{j+1/2}}{\|[[\nabla_{\bf u}\eta_\btheta({\bf u})]]_{j+1/2}\|^2_2} [[\nabla_{\bf u}\eta_\btheta({\bf u})]]_{j+1/2}\\
&\ - \frac{\frac{1}{2}[[\nabla_{\bf u}\eta_\btheta({\bf u})]]^\top\big({\bf f}^{\btheta,\bmu}({\bf u}_{j+1/2}^+) + {\bf f}^{\btheta,\bmu}({\bf u}_{j+1/2}^-)\big)}{\|[[\nabla_{\bf u}\eta_\btheta({\bf u})]]_{j+1/2}\|^2_2} [[\nabla_{\bf u}\eta_\btheta({\bf u})]]_{j+1/2},
\end{aligned}
\end{equation}
where we have substituted $({\bf f},\phi,{\bf v})$ in \eqref{eq: entropy conservative flux} by $({\bf f}^{\btheta,\bmu},\phi_\bmu,\nabla_{\bf u}\eta_\btheta({\bf u})))$. Here, ${\bf f}^{\btheta,\bmu}$ is defined in \eqref{eq: paramatric HCLaw}, and $\mathbf{u}_{j+1/2}^\pm$ are reconstructed interface states. For completeness we detail the construction of the left-biased state ${\bf u}_{j+1/2}^-$ for the WENO5 scheme \cite{shu2020essentially}, noting that the right-biased state ${\bf u}_{j+1/2}^+$ is obtained by symmetry. Specifically, ${\bf u}_{j+1/2}^-$ is expressed as a convex combination of three third-order candidate polynomials with
\begin{align}
\bu_{j+\frac{1}{2}}^- = w_0 p_0 + w_1 p_1 + w_2 p_2,
\label{eq: weno5_l}
\end{align}
Here the third-order candidate polynomials are given by
\[\begin{aligned} 
p_0 &= \tfrac{1}{3}\bu_{j-2} - \tfrac{7}{6}\bu_{j-1} + \tfrac{11}{6}\bu_j,\\
p_1 &= -\tfrac{1}{6}\bu_{j-1} + \tfrac{5}{6}\bu_j + \tfrac{1}{3}\bu_{j+1}, \\
p_2 &= \tfrac{1}{3}\bu_j + \tfrac{5}{6}\bu_{j+1} - \tfrac{1}{6}\bu_{j+2},\end{aligned}\]
and the nonlinear weights are defined as $w_k = \alpha_k/(\alpha_1+\alpha_2+\alpha_3)$ with $\alpha_k = c_k/((\Delta x)^2+\beta_k)^2$, $k = 1,2,3$.   The corresponding smoothness indicators are determined as
\[\begin{aligned}
\beta_0 &= \tfrac{13}{12}(\bu_{j-2} - 2\bu_{j-1} + \bu_j)^2 + \tfrac{1}{4}(\bu_{j-2} - 4\bu_{j-1} + 3\bu_j)^2, \\
\beta_1 &= \tfrac{13}{12}(\bu_{j-1} - 2\bu_j + \bu_{j+1})^2 + \tfrac{1}{4}(\bu_{j-1} - \bu_{j+1})^2, \\
\beta_2 &= \tfrac{13}{12}(\bu_j - 2\bu_{j+1} + \bu_{j+2})^2 + \tfrac{1}{4}(3\bu_j - 4\bu_{j+1} + \bu_{j+2})^2,
\end{aligned}\]
and  the optimal weights are prescribed as $c_0 = 0.1$, $c_1 = 0.6$ and $c_2 = 0.3$.

Combining \eqref{eq: Param entropy stable flux}, \eqref{eq:entropy_stable_flux_nn}--\eqref{eq: weno5_l}, we obtain an entropy-stable flux for SymCLaw, 
\begin{equation}
   \widehat{\bf f}_{j+1/2}^{\btheta,\bmu} 
   =  {\bf f}^{\btheta,\bmu, \ast}_{j+1/2} \!\!- \tfrac{1}{2}  \lambda_{j+1/2}^\text{max}\big(\mathbb{H}_{\bf u}\eta_{\btheta}(\bar {\bf u}_{j+1/2})\big)^{-1}\!\!\left(\!\nabla_{\bf u}\eta_\btheta({\bf u}_{j+1/2}^+) \!-\! \nabla_{\bf u}\eta_\btheta({\bf u}_{j+1/2}^-)\!\right),
\label{eq: neural-entropy stable flux}
\end{equation}
where $\lambda_{j+1/2}^{\text{max}}$ is defined in \eqref{eq:maximum_wave_speed} and depends on $\bmu$ and $\btheta$. The details of the structure of neural networks and the choice of training parameters are deferred to \Cref{sec:networkdetails}.

 With the entropy-stable flux \eqref{eq: neural-entropy stable flux} in hand, the semi-discretized problem formulation for the SymCLaw \eqref{eq: paramatric HCLaw} is given by
\begin{equation}\label{eq: scheme_learned_law}
\frac{d}{dt} {\bf u}_j(t) = -\frac{1}{\Delta x} \left( \widehat{\bf f}_{j+1/2}^{\btheta,\bmu} - \widehat{\bf f}_{j-1/2}^{\btheta,\bmu}\right).
\end{equation}
Thus we see that identifying the parameters $(\btheta, \bmu)$ reduces to minimizing the discrepancy between observed data trajectories and those generated by the model \eqref{eq: scheme_learned_law}. 

\subsection{Time integration} \label{sec: time_integrator} To ensure consistency with the spatial discretization order, where WENO5 reconstruction is employed for the interface states, we adopt the third-order Total Variation Diminishing Runge--Kutta (TVDRK3) scheme (\cite{Shu98}) for time integration of the  system of ODEs obtained in \cref{eq: scheme_learned_law}.  This method can be applied to any generic time-dependent ODE of the form
$$ \frac{dz}{dt}= \mathcal{L}(z),$$
where $\mathcal{L}$ is a known operator acting on $z$. The TVDRK3 algorithm, summarized in \Cref{alg: TVDRK3}, advances the solution $z^{l-1}$ at time ${t}_{l-1}$ to $z^l$ at the next time level ${t}_l$, for $l \geq 1$. The time step $\Delta t$ is fixed in all of our numerical experiments for simplicity.

\begin{algorithm}[h!]
\caption{TVDRK3 time integration method for a single time step starting at time level $t^{l-1}$}\label{alg: TVDRK3}
\begin{algorithmic}
\State INPUT: ${z}^{l-1}$, $\mathcal{L}(z^{l-1})$ and $\Delta t$
\State OUTPUT: The solution ${z}^{l}$ at time level $t_l$
\State $ z^{(1)}=z^{l-1}+ \Delta t \, \mathcal{L}(z^{l-1})$
\State $ z^{(2)}=\frac{3}{4}z^{l-1} + \frac{1}{4} z^{(1)} + \frac{1}{4} \Delta t \, \mathcal{L}(z^{(1)})$
\medskip
\State $ z^{l}=\frac{1}{3}z^{l-1} + \frac{2}{3}z^{(2)} + \frac{2}{3}\Delta t \, \mathcal{L}(z^{(2)})$
\end{algorithmic}
\end{algorithm}

\section{Training procedure design}\label{sec: training design}

 We now provide details for the data generation process and  training protocol used in our numerical experiments. 

\subsection{Data generation}\label{sec: data generation} 
\begin{figure}[ht]
\centering
\includegraphics[width=0.7\textwidth,trim={0 0 0 0},clip]{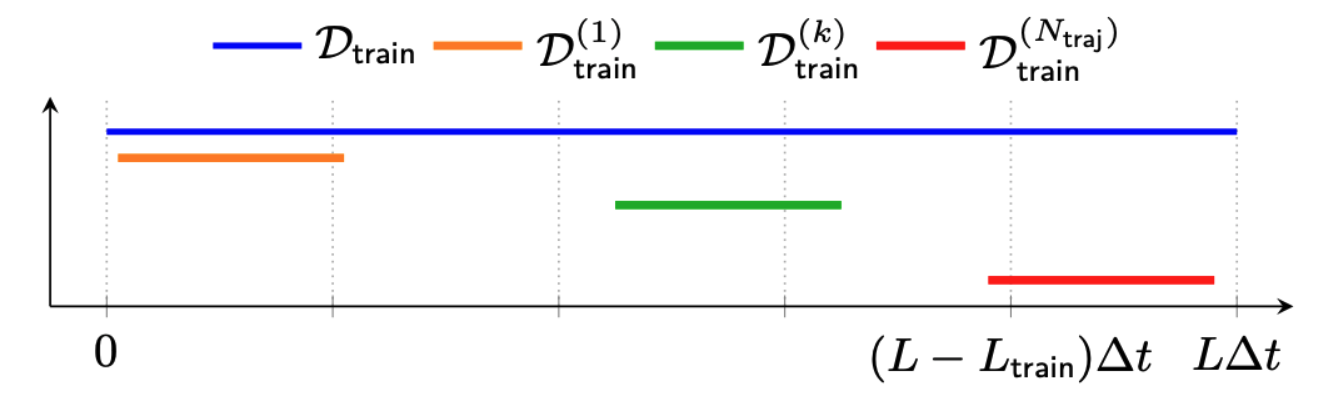}
\caption{Illustration of full time span and training subintervals.}\label{fig: training domain}
\end{figure}

We assume access to training data in the form of solution trajectories generated over a finite time interval. Each trajectory originates from a perturbed initial state and evolves according to the governing conservation law. These trajectories are discretized in both time and space. In practice, such data may arise from experimental measurements or sensor networks. In this study, however, we simulate observations by numerically solving the true PDE using various perturbed initial conditions. Importantly, and in contrast to some prior work (e.g., \cite{chenPDE}), we do {\em not} rely on oracle access that might provide a richer solution space but would not be available as observations, nor do we carefully curate initial conditions to guarantee smoothness of the solution. Instead, we adopt a more realistic and general data acquisition framework, in line with the assumptions in \cite{chencfn,liu2024entropy}. Specifically, we discretize the temporal domain using a fixed time step $\Delta t$ and define by $L$ the total number of simulation steps. The full time span for data collection is then given by
\begin{equation}
    \label{eq:totalperiod}
    \mathcal{D}_{\text{train}}= [0,L\Delta t],
\end{equation} 
over which we generate trajectories from $N_{\text{traj}}$ different initial conditions. For each trajectory indexed by $k \in \{1,\cdots, N_{\text{traj}}\}$, we extract a training subinterval
\begin{equation}
    \label{eq:trainingperiod}
    \mathcal{D}^{(k)}_{\text{train}}= [t_0^{(k)},t^{(k)}_{L_{\text{train}}}],
\end{equation} 
where the starting time $t_0^{(k)}$ is sampled from the interval $[0, (L-L_{\text{train}})\Delta t]$, and the terminal time is given by $t^{(k)}_{L_{\text{train}}} = t_0^{(k)} + L_{\text{train}}\Delta t$, see \Cref{fig: training domain} for an illustration. 

The observed data along each trajectory then consists of discrete space-time samples of the conserved states
\begin{equation}
\label{eq: trajectory data}
\bu(t_l^{(k)}) \in {\mathbb R}^{n_{\text{train}}\times p}, \ l=0, \cdots, L_{\text{train}}, \ k=1, \cdots, N_{\text{traj}},
\end{equation}
where $n_{\text{train}}$ denotes the number of spatial grid points used during training, and $p$ is the number of the conserved states.

\subsection{Loss function} \label{sec:recurrentloss} 
For ease of presentation, we denote by $\mathcal{N}$ the procedure that combines \eqref{eq: scheme_learned_law} with the TVDRK3 time integrator in \Cref{alg: TVDRK3}, representing the neural operator for advancing the solution in time. That is, given the current prediction at $t_{l}^{(k)}$, $\hat{\bu}(t_l^{(k)}) \in \mathbb{R}^{n_{\text{train}} \times p}$, the next-step prediction is
\begin{equation}
    \label{eq:Ndefine}
    \hat\bu(t_{l+1}^{(k)}) = \mathcal{N}\big(\hat\bu(t_{l}^{(k)})\big) \quad \text{ with } \quad \hat\bu(t_0^{(k)}) = \bu(t_0^{(k)}).
\end{equation}
 We can now define the normalized \emph{recurrent loss function} as
\begin{equation}
\label{eq: stage 1 loss}
\mathfrak{L}\left( \btheta, \bmu; {\bf u} \right) = \frac{\sum_{k=1}^{N_{\text{traj}}} \sum_{l=0}^{L_{\text{train}}} \left\|\hat\bu( t_{l}^{(k)} ;\btheta, \bmu) - \bu( t_{l}^{(k)} )\right\|_1}{\sum_{k=1}^{N_{\text{traj}}} \sum_{l=0}^{L_{\text{train}}} \|\bu( t_{l}^{(k)} )\|_1}, 
\end{equation}
 where 
 \[\hat\bu( t_{l}^{(k)};\btheta, \bmu) = \underbrace{\mathcal{N}\circ\ldots \circ \mathcal{N}}_{l \text{ times }}\big( \bu( t_{0}^{(k)}) \big).\]
\begin{remark}
In contrast to the NESCFN approach in \cite{liu2025neural}, which enforces structural properties such as symmetry through regularization terms in the loss function and therefore achieves hyperbolicity only in an empirical sense, the newly proposed SymCLaw method embeds these properties directly into the network architecture. Consequently, the loss function involves only data misfit terms, which simplifies the training procedure and eliminates dependence on regularization parameters. Moreover, SymCLaw admits a rigorous entropy stability guarantee: an entropy-conservative flux \eqref{eq:entropy_stable_flux_nn} is constructed by design and subsequently augmented with appropriate numerical dissipation \eqref{eq: Param entropy stable flux}. Such a theoretical entropy stability guarantee is not available for NESCFN, where entropy stability is observed numerically but not ensured by construction.
\end{remark}

\subsection{Conservation and entropy metric}\label{sec:conservationnmetric}
Motivated by the evaluation strategy in \cite{liu2025neural}, we examine whether SymCLaw preserves the conservation of physical quantities and maintains entropy stability throughout the simulation. 

Precisely, to quantify conservation fidelity, we define the discrete conservation error for the $\mathfrak i$-{th} conserved variable at a future time $t_l$ beyond the training horizon as
\begin{equation}
\label{eq:conservemetric}
\mathcal{C}_{\mathfrak{i}}\big(\hat\bu(t_l)\big) := \Big| \sum_{j=1}^{n_{\text{train}}-1} \left( \hat{u}_j^{\mathfrak{i}}(t_l) - \hat{u}_j^{\mathfrak{i}}(t_0) \right)\Delta x - \sum_{s=1}^{l}\left( F_a^{{\mathfrak{i}},s-1} -F_b^{{\mathfrak{i}},s-1} \right)\Delta t\Big|,
\end{equation}
where $[\hat u^{\mathfrak{i}}_1(t_l), \ldots, \hat u^{\mathfrak{i}}_{n_{\text{train}}-1}(t_l)]$ denotes the predicted ${\mathfrak{i}}$-th conserved state at time \(t_l\), and \(F_a^{{\mathfrak{i}},s-1}\) and \(F_b^{{\mathfrak{i}},s-1}\) are the calculated fluxes at the domain boundaries $x = a$ and $x = b$, respectively, defined as
\begin{equation}
\label{eq:fluxterm}
        F_a^{\mathfrak{i},s-1} = \frac{1}{\Delta t}\int_{t_{s-1}}^{t_{s}} {\bf f}^{\bmu,\btheta}_\mathfrak{i}\big(\hat{\bf u}(a,t)\big)dt,\quad
        F_b^{\mathfrak{i},s-1} = \frac{1}{\Delta t}\int_{t_{s-1}}^{t_{s}} {\bf f}_{\mathfrak{i}}^{\bmu,\btheta}\big(\hat{\bf u}(b,t)\big)dt.
\end{equation}
Here ${\bf f}_{\mathfrak{i}}^{\bmu,\btheta}$ denotes the $\mathfrak{i}$-th component of the learned flux defined in \eqref{eq: paramatric HCLaw} and $n_{\text{train}} = n_{\text{test}}$, the number of spatial grid points used for testing in \Cref{sec: numerics}.

To assess the entropy stability, we consider the discrete entropy remainder at $t_l$ 
\begin{equation}
\begin{aligned}
\mathcal{J}(\hat\bu (t_l)) 
:=   \sum_{j=1}^{n_{\text{train}}-1}\bigg[ &\Delta x \Big( {\eta}_{\btheta}\big(\hat{\bf u}_j(t_l)\big) - {\eta}_{\btheta}\big(\hat{\bf u}_j(t_0)\big) \Big) \\
&- \sum_{s=1}^{l} \left( \big(\nabla_{\bf u} \eta_{\btheta}\big(\hat{\bf u}_{j}(t_s)\big)\big)^\top{\bf f}^{\bmu,\btheta}\big( \hat{\bf u}_{j}(t_s)\big) \right) \Delta t \bigg],
\end{aligned}
    \label{eq: discrete entropy}
\end{equation}
where the second term approximates the entropy flux. We say that the learned conservation law is entropy stable if $\mathcal{J}(\hat\bu (t_l)) \leq 0$. 

\subsection{Network architectures and training configuration}
\label{sec:networkdetails}

With the data preparation and training metrics established in the previous sections, we now describe the neural network architectures employed in SymCLaw, together with the associated training parameters and stabilization techniques. 

As discussed in \Cref{sec: parametric conservation laws}, the SymCLaw \eqref{eq: paramatric HCLaw} is realized by matching observed experimental data with predictions generated from the corresponding discrete model. This formulation involves two types of neural networks. Specifically, the entropy function $\eta_{\btheta}$ is parameterized by an input convex neural network, ensuring the convexity of the learned entropy. The flux function, expressed in terms of the entropy variable, is implemented as the \textit{gradient} of a fully connected neural network, namely, ${\bf g}_{\mathfrak{i}}({\bf v}) = \nabla_{\bf v} \phi_{\bmu,\mathfrak{i}}({\bf v})$. This construction guarantees that its Jacobian is symmetric, which, together with the convexity of $\eta_{\btheta}$, ensures the hyperbolicity of the SymCLaw. For completeness, we summarize the architectures of these networks below.

\subsubsection{Fully connected neural networks} A fully connected neural network (FCNN) is a feedforward architecture composed of multiple layers, where each layer applies an affine transformation followed by a nonlinear activation. Let ${\bf z}^{(0)} = \mathbf{x} \in \mathbb{R}^{d_0}$ denote the input. An $\mathrm{L}$-layer FCNN \(\phi_\bmu(\mathbf{x})\) computes the output $\mathbf{y} \in \mathbb{R}^{\mathrm{L}}$ through the recursive relation
\begin{equation}
\begin{aligned}
\mathbf{z}^{(l)} &= \sigma\Big(\mathbf{W}^{(l)} \mathbf{z}^{(l-1)} + \mathbf{b}^{(l)}\Big), \quad l = 1,\cdots, \mathrm{L}-1,\\
{\bf y} & = \mathbf{W}^{(\mathrm{L})} \mathbf{z}^{(\mathrm{L}-1)} + {\bf b}^{(\mathrm{L})},
\end{aligned}\label{eq:activation}
\end{equation}
where $\mathbf{z}^{(0)} = \mathbf{x}$ is the input, $\mathbf{W}^{(l)} \in \mathbb{R}^{d_l \times d_{l-1}}$ and $\mathbf{b}^{(l)} \in \mathbb{R}^{d_l}$ are the weight matrices and bias vectors at layer $l$, and $\sigma(\cdot)$ is an elementwise nonlinear activation function. 

Our SymCLaw method employs the gradients of FCNNs  with respect to the input \eqref{eq:activation} to parameterize the flux function $\mathbf{f}(\mathbf{u})$ in terms of the entropy variable $\mathbf{v}$, see \eqref{eq:gradient_NN}.

\subsubsection{Input convex neural networks}\label{sec:input_cnn}  
Input convex neural networks (ICNNs) form a subclass of FCNNs designed to represent functions that are convex with respect to their input. Introduced in~\cite{amos2017input}, ICNNs impose two key structural constraints:  
\begin{itemize}
\item[i.] The activation functions must be convex and non-decreasing;  
\item[ii.] The weight matrices in the hidden-layer recursion must be element-wise non-negative to preserve convexity.  
\end{itemize}  

Formally, for an ICNN with input ${\bf z}^{(0)} = \mathbf{x} \in \mathbb{R}^{d_0}$ and layer outputs $\mathbf{z}^{(l)} \in \mathbb{R}^{d_l}$, the forward pass is given by  
\begin{equation} 
\begin{aligned}
{\bf z}^{(l)} &= \sigma\!\left({\bf W}^{z,(l)} {\bf z}^{(l-1)} + {\bf W}^{x,(l)} {\bf x} + {\bf b}^{(l)}\right), \quad l = 1,\ldots,\mathrm{L}-1,\\
{\bf y} &= \mathbf{W}^{(\mathrm{L})} \mathbf{z}^{(\mathrm{L}-1)} + {\bf b}^{(\mathrm{L})},
\end{aligned}\label{eq:ICNNstructure}
\end{equation}
where ${\bf y} \in \mathbb{R}^{d_\mathrm{L}}$ is the output. To enforce convexity, each hidden-layer weight matrix $\mathbf{W}^{z,(l)}$ is projected elementwise onto the nonnegative orthant, i.e.,
\[\mathbf{W}^{z,(l)} \leftarrow \max(\mathbf{W}^{z,(l)},0),\quad  l = 1,\ldots,\mathrm{L}-1.\] 
We further augment the output layer with quadratic and linear terms in the input
\[{\bf y} = \mathbf{W}^{(\mathrm{L})} \mathbf{z}^{(\mathrm{L}-1)} + {\bf b}^{(\mathrm{L})} + \mathbf{W}_{\mathrm{S}}\mathbf{x}^2 + \mathbf{W}_{\mathrm{L}}\mathbf{x}.\]
The addition of $\mathbf{W}_{\mathrm{S}}\mathbf{x}^2 + \mathbf{W}_{\mathrm{L}}\mathbf{x}$ increases the expressive power of the ICNN by enabling it to directly represent quadratic and linear contributions. 
To ensure convexity, the quadratic weight matrix $\mathbf{W}_{\mathrm{S}}$ is regularized using the Huber penalty:
\[
 \mathbf{W}_{\mathrm{S}} =
 \begin{cases}
 |\mathbf{W}_{\mathrm{S}}| - 0.5, & |\mathbf{W}_{\mathrm{S}}| > 1,\\[6pt]
 \tfrac{1}{2}\mathbf{W}_{\mathrm{S}}^2, & |\mathbf{W}_{\mathrm{S}}| \leq 1.
 \end{cases}
\]  

ICNNs are particularly valuable when convexity of the learned function is a modeling requirement. In our framework, this property is crucial because the entropy function $\eta_{\btheta}(\mathbf{u})$ must be convex in order to guarantee the hyperbolicity of the learned conservation law. 
To this end, we explicitly employ an ICNN of the form \eqref{eq:ICNNstructure} to parametrize and learn the entropy function $\eta_\btheta(\bu)$, where $\btheta$ collectively denotes all trainable parameters in \eqref{eq:ICNNstructure}.

 \subsubsection{Training configuration} We now describe the training configuration used in this work for the numerical experiments conducted in \Cref{sec: numerics}. As discussed in the preceding sections, the neural operator responsible for solution updates, denoted by $\mathcal{N}$ in~\eqref{eq:Ndefine}, is implemented using fully connected neural networks. Each network has both input and output dimensions equal to the number of state variables $\bf u$, i.e., $d_0 = d_{\mathrm{L}} = p$ in \eqref{eq:activation}--\eqref{eq:ICNNstructure}. The architectural details of the individual networks for different numerical examples in \Cref{sec: numerics} are summarized as follows:
 \begin{enumerate}
 \item [i.] The FCNN $\phi_{\bmu,\mathfrak{i}}$ consists of $3$ hidden layers, each with $64$ neurons for both the 1D shallow water and 1D Euler's equations; and $1$ hidden layer with $32$ neurons for the scalar hyperbolic conservation laws, which include both the 1D and 2D Burgers' equation, as well as the 2D KPP equation. Such architectures balance expressiveness with memory constraints. We employ the tanh activation function, defined as $\text{Tanh}(x) = \tanh(x)$ in the hidden layers. The choice of tanh is motivated by its differentiability, which is essential for computing the Hessian of $\phi_{\bmu,\mathfrak{i}}$ in \eqref{eq:maximum_wave_speed}.
\item [ii.] The ICNN $\eta_\btheta$ is modeled using a two-layer ICNN with $64$ hidden neurons each layer for systems of conservation laws (1D shallow water  and 1D Euler's equations); and $2$ hidden layer with $32$ neurons for scalar hyperbolic conservation laws (1D and 2D Burgers' equation, and 2D KPP equation). To ensure convexity and smoothness, we apply the Softplus activation function, defined as $\text{Softplus}(x) = \log(1+e^x)$, which is convex, non-decreasing and continuously differentiable to ensure the convexity of $\eta_\btheta$. 
 \end{enumerate}

To determine the parameters $(\btheta,\bmu)$ in \cref{eq: paramatric HCLaw}, we employ the Adam optimizer (\cite{Adam}) to update the neural networks' weights and biases during training. We adopt a one-cycle cosine learning-rate scheduler throughout to ensure consistency across experiments. Specifically, the peak learning rate is set to \(5\times10^{-3}\); the initial and final division factors are chosen as \(10\) and \(10^3\), respectively; and the warm-up phase is taken as the first \(10\%\) of the training steps.\footnote{Because of its greater complexity, we extend the warm-up phase to the first \(5\%\) of the training steps for the Euler equations.}   \Cref{tab:training parameter} provides the set of all other training hyperparameters for experiments conducted in \Cref{sec: numerics}: the number of training epochs $N_{\text{Epoch}}$; the number of training trajectories \(N_{\text{traj}}\); and the batch size $N_b$. These choices reflect a balance between problem complexity, training efficiency, and available GPU memory. While the number of training trajectories, \(N_{\text{traj}}\), varies across problems, we fix the number of validation cases to $40$ in all experiments, which is used to determine whether the model parameters should be saved at each epoch. 
\begin{table}[h!]
\centering
\scalebox{0.92}{
\begin{tabular}{|l|c|c|c|c|c|c|}
\hline
\text{numerical examples} & \small{$N_{\text{Epoch}}$} & \(N_{\text{traj}}\) & $N_b$\\
\hline
\small{1D Burgers' eq}   &  $200$ & $200$ & 5\\
\hline
\small{1D shallow water eq}   & $200$ & $300$ & 10\\
\hline
\small{1D Euler's eq}    & $500$ & $150$  & 5\\
\hline
\small{2D Burgers' eq}   & $500$ & $10$ & 2\\
\hline
\small{2D KPP eq}   & $500$ & $50$ & 2\\
\hline
\end{tabular}}
\caption{Training parameters for each experiment conducted in \Cref{sec: numerics}.}
\label{tab:training parameter}
\end{table}

Finally, we note that all implementations are based on JAX \cite{jax}, a high-performance numerical computing library in Python that supports automatic differentiation and GPU/TPU acceleration and that all experiments are implemented in double precision. To ensure robustness and fairness, none of the training parameters in our method were further tuned.\footnote{The complete code is available upon request for reproducibility.}

\subsubsection{Techniques to stabilize training} Our training procedure discussion concludes by outlining the stabilization techniques employed to ensure robustness of the learning process. Constructing the entropy-stable flux in \eqref{eq: Param entropy stable flux} requires both computing the maximal local wave speed $\lambda_{j+1/2}^{\text{max}}$ and inverting the entropy Hessian $\mathbb{H}_{\bf u}\eta_{\btheta}$. However, two challenges arise during the early stages of training: (i) the Hessian produced by the neural network can be highly unstable, either exhibiting excessively large eigenvalues or degenerating with entries close to zero; and (ii) the fixed time-step size $\Delta t$ may violate the CFL condition when combined with the neural-network-predicted local wave speed. To address these issues, we employ three stabilization strategies:
\begin{enumerate}
\item [i.] {\bf Epoch-dependent Hessian regularization.} We add a diagonal perturbation with exponentially decaying magnitude to improve the conditioning of the Hessian:
\begin{equation*}
\widetilde{\mathbb{H}}_{\bf u}\eta_{\btheta}(\bar{\bf u}_{j+1/2}) = \mathbb{H}_{\bf u}\eta_{\btheta}(\bar{\bf u}_{j+1/2}) + C_1^{\text{Epoch} -1}\mathbb{I}_p, 
\label{eq: stable cho decomp}
\end{equation*}
where $\mathbb{I}_p$ is the $p\times p$ identity matrix, and $0 < C_1 < 1$. We choose $C_1 = 0.1$ to accelerate the warm-up process of neural networks.

\item [ii.] {\bf Stabilized flux correction.} Invoking the relation between entropy and state variables, we have 
    \[\big(\mathbb{H}_{\bf u}\eta(\bar{\bf u}_{j+1/2})\big)^{-1}[[{\bf v}]]_{j+1/2} = [[{\bf u}]]_{j+1/2}.\]
    However, since $\widetilde{\mathbb{H}}_{\bf u}\eta_{\btheta}$ is only an approximation, we enforce a correction rule:
    \[
    \big(\widetilde{\mathbb{H}}_{\bf u}\eta_{\btheta}(\bar{\bf u})\big)^{-1}[[\bf v]] = \left\{ 
    \begin{aligned}
        &\big(\widetilde{\mathbb{H}}_{\bf u}\eta_{\btheta}(\bar{\bf u})\big)^{-1}[[\bf v]],  & \| \big(\widetilde{\mathbb{H}}_{\bf u}\eta_{\btheta}(\bar{\bf u})\big)^{-1}[[{\bf v}]]\|_{\infty}\leq C_{D}\|[[{\bf u}]]\|_{\infty}\\
        & [[{\bf u}]], \quad& \text{otherwise}
    \end{aligned},
    \right.
    \]
   where $C_D > 1$ (here we choose $C_D= 2$) is a fixed constant.  The subscript $j+1/2$ is omitted for clarity. We note that the worst-case scenario for this  stabilization strategy results in the reduction of \eqref{eq: Param entropy stable flux} to the local Lax--Friedrichs flux. 
    
\item [iii.] {\bf Wave-speed clipping.} To avoid violating the CFL condition, we clip the estimated local maximum wave speed using prior knowledge of the observed data:
    \[\tilde{\lambda}^{\text{max}}_{j+1/2} = \min\left(\lambda_{j+1/2}^{\text{max}}, C_{\text{CFL}}\frac{\Delta x}{\Delta t}\right).
    \]  
We emphasize that this CFL safeguard is merely a workaround. A more principled approach would be to adaptively determine the time-step size, which would require differentiable implementations of adaptive solvers for backpropagation. We leave such extensions to future work. In our experiments we choose $C_{\text{CFL}} = 1$.
\end{enumerate}

\section{Numerical results}\label{sec: numerics}

We now present numerical experiments to demonstrate the effectiveness of our new SymCLaw methodology. In particular, we show that our approach not only predicts the long-term behavior of the dynamics with high fidelity, but also preserves entropy stability across all levels of noise  present in the training data. 

\subsection{Prototype conservation laws}
\label{sec: pdeexamples}
The classical conservation laws introduced below serve as benchmark test cases for evaluating the performance of the SymCLaw.

\subsubsection{1D Burgers' equation}
\label{sub:1DBurgers}
We first consider the scalar Burgers' equation
\begin{equation}
    u_t + \Big( \dfrac{u^2}{2} \Big)_x = 0, \quad x \in \left[ 0, 2\pi \right], \quad t > 0,
    \label{eq: Burgers' equation}
\end{equation}
with periodic boundary conditions, and initial condition
\begin{equation}
    u\left( x, 0 \right) = \alpha \sin\left( x \right) + \beta, \quad \alpha,\beta \in \mathbb{R}.
    \label{eq: Burgers' initial condition}
\end{equation}

We employ the PyClaw package \cite{clawpack, pyclaw} with a fixed time step \(\Delta t = .005\) and spatial resolution $n_{\text{train}} = 512$ for training.  The initial data \eqref{eq: Burgers' initial condition} are generated by sampling $\alpha \sim \left[ .75, 1.25 \right]$ and $\beta \sim \left[ -.25, .25 \right]$ uniformly. The total training time horizon is set to $L = 20$ steps (see \cref{eq:totalperiod}), corresponding to $\mathcal{D}_{\text{train}} = [0,L\Delta t] = [0,0.1]$. In this example we do not subdivide the training interval, i.e., $L_{\text{train}} = L =20$ in \cref{eq:trainingperiod}. Moreover, we emphasize that only \emph{smooth} solution profiles prior to shock formation are available for training.

For testing, we fix $\alpha = 1.05609$ and $\beta = 0.1997$ in \cref{eq: Burgers' initial condition}, and compute the reference solution up to $T = 3$ using $(\Delta t, n_{\text{test}}) = (.005, 512)$.

\subsubsection{Shallow water system} \label{sub:shallowwater} We next consider 1D shallow water system
\begin{equation}
    \begin{aligned}
        h_t + \left( hu \right)_x &= 0, \\
        \left( hu \right)_t + \Big( hu^2 + \dfrac{1}{2}gh^2 \Big)_x &= 0,
    \end{aligned}
    \label{eq: shallow water equation}
\end{equation}
over $x\in(-5,5)$ with Dirichlet boundary conditions. The initial conditions follow
\begin{equation}
    h\left( x, 0 \right) = \left\{
    \begin{aligned}
     &h_l + \omega_{l}, \quad x < x_0 + \omega_{x}   \\
     &h_r + \omega_{r}, \quad x \geq x_0 + \omega_{x} \\
     \end{aligned}
    \right.,\quad 
    u\left( x, 0 \right) = \left\{
    \begin{aligned}
     &u_l+\omega_{ul}, \quad x < x_0+\omega_{x}   \\
     &u_r+\omega_{ur}, \quad x \geq x_0+\omega_{x} \\
     \end{aligned}
    \right.,\label{eq: shallow water condition}
\end{equation}
where \(h_l = 3.5, h_r=1.0, u_l=u_r=x_0=0\), and \(\omega_{l}, \omega_{r}, \omega_{ul}, \omega_{ur}, \omega_{x} \in \mathbb{R} \).

The initial conditions for the training trajectories are given by \eqref{eq: shallow water condition}, where \(\omega_{l}, \omega_{r}\sim \mathcal{U}\left[ -.2, .2 \right]\) and \(\omega_{ul}, \omega_{ur}, \omega_{x} \sim \mathcal{U}\left[ -.1,.1 \right]\). All other training parameters, \((\Delta t, n_{\text{train}}, L, L_{\text{train}})\), are equivalent to those used for the 1D Burgers' equation in \Cref{sub:1DBurgers}.  

For testing, we set  
$\omega_{l}=\omega_{r}=\omega_{ul}=\omega_{ur}=\omega_{x}=0$
and fix
\[(h_l,h_r,u_l,u_r,x_0) = (3.5691196,1.178673, -.064667, -.045197, 003832)\] 
in \eqref{eq: shallow water condition}. The reference solution is computed with the same temporal and spatial resolution as in the training data and is simulated up to $T = 1.5$.

\subsubsection{Euler's equations}\label{sub:Eulers} We now consider the system of 1D Euler's equations\begin{equation}
    \begin{aligned}
        \rho_{t} + \left( \rho u \right)_{x} &= 0, \\
        \left( \rho u \right)_{t} + \Big( \rho u^2 + p \Big)_{x} &= 0, \\
        \left( E \right)_{t} + \left( u (E + p) \right)_{x} &= 0,
    \end{aligned}
    \label{eq: Euler's equation}
\end{equation}
over the spatial domain $x \in (-5,5)$ with Dirichlet boundary conditions. We assess the performance of SymCLaw on the challenging Shu--Osher problem with initial conditions
\begin{equation}
\begin{aligned}
& \rho(x, 0)= \begin{cases}\rho_l, & \text { if } x \leq x_0, \\
1+\varepsilon \sin (5 x), & \text { if } x_0<x \leq x_1, \quad u(x, 0)= \begin{cases}u_l, & \text { if } x \leq x_0, \\
0, & \text { otherwise, }\end{cases} \\
1+\varepsilon \sin (5 x) e^{-\left(x-x_1\right)^4}, & \text { otherwise, }\end{cases} \\ \\
& p(x, 0)=\left\{
\begin{array}{ll}
p_l, & \text { if } x \leq x_0, \\
p_r, & \text { otherwise, }
\end{array} \quad E(x, 0)=\frac{p_0}{\gamma-1}+\frac{1}{2} \rho(x, 0) u(x, 0)^2, \right.
\end{aligned}\label{eq:euler IC}
\end{equation}
where $x_1=3.29867$ and $\gamma=1.4$. 

To generate training trajectories, we uniformly sample parameters in \eqref{eq:euler IC} as\vspace{-0.2cm}
\[\begin{+array}{lrlr}
 \rho_l \sim \mathcal{U}[\hat{\rho}_l(1-\epsilon), &\hat{\rho}_l(1+\epsilon)], &
 \varepsilon\; \sim \mathcal{U}[\hat{\varepsilon}(1-\epsilon),  &\hat{\varepsilon}(1+\epsilon)], \\
 p_l \sim \mathcal{U}[\hat{p}_l(1-\epsilon), &\hat{p}_l(1+\epsilon)], &
 p_r \sim \mathcal{U}[\hat{p}_r(1-\epsilon),  &\hat{p}_r(1+\epsilon)], \\
 u_l \sim \mathcal{U}[\hat{u}_l(1-\epsilon), &\hat{u}_l(1+\epsilon)], & 
 x_0 \sim \mathcal{U}[\hat{x}_0(1-\epsilon), &\hat{x}_0(1+\epsilon)],
\end{+array}\vspace{-0.2cm}\]
where $(\epsilon, \hat\rho_l, \hat p_l, \hat u_l, \hat \varepsilon, \hat p_r, \hat x_0)=(.1,3.857135,10.32333,2.62936,.2,1,-4)$. Each trajectory is obtained by solving \eqref{eq: Euler's equation} with PyClaw using the HLLE Riemann solver and spatial-temporal resolution \((\Delta t,n_{\text{train}}) = (.002, 512)\), evolved until $T = L\Delta t$. 

Due to the complexity of the Euler's equations and limitations on GPU memory,  here we  subdivide the training interval $\mathcal{D}_{\text{train}} = [0, L\Delta t]$ into $N_{\text{traj}}$ subintervals (see \Cref{fig: training domain}).  Specifically, each $k$-th trajectory is associated with training interval $\mathcal{D}_{\text{train}}^{(k)} = [t_0^{(k)}, t_0^{(k)} + L_{\text{train}}\Delta t]$. The  starting time is given by $t_0^{(k)}\sim \mathcal{U}\{0,1, \cdots, L-L_{\text{train}}\}$, and in this example we choose $L = 300$ and $L_{\text{train}} = 20$.

For testing, we set \((\rho_l,\varepsilon,p_l,p_r,u_l,x_0) = (\hat{\rho}_l, \hat{\varepsilon}, \hat{p}_l, \hat{p}_r, \hat{u}_l, \hat{x}_0)\) in \eqref{eq:euler IC}, and simulate the system \eqref{eq: Euler's equation} until \(T = 1.6\) with $(\Delta t,n_{\text{test}}) = (.002, 512)$.

\subsubsection{2D Burgers' equation}\label{sub: 2D Burgers} To evaluate the SymCLaw method in multidimensional settings, we first consider the 2D Burgers' equation
\begin{equation}
    u_t + \Big( \dfrac{u^2}{2} \Big)_x + \Big( \dfrac{u^2}{2} \Big)_y = 0, \quad (x, y) \in \left[ 0, 1 \right]\times\left[ 0, 1 \right], \quad t > 0, 
    \label{eq: Burgers' equation 2d}
\end{equation}
subject to periodic boundary conditions and initial condition  given by
\begin{equation}   
u\left( x, y, 0 \right) = \alpha \sin\left( 2\pi x +x_0\right)\cos\left( 2\pi y +y_0\right) + \beta,\quad \alpha, \beta,x_0,y_0 \in \mathbb{R}.
    \label{eq: Burgers' initial condition 2d}
\end{equation}

For training, the parameters in \eqref{eq: Burgers' initial condition 2d} are sampled uniformly with \(\alpha \sim {\mathcal U}[.75,1.25]\), \(\beta \sim {\mathcal U}[-.25,.25]\),  \(x_0 \sim {\mathcal U}[.5,1.5]\), and \(y_0 \sim {\mathcal U}[-.5, .5].\) Each of the $N_{\text{traj}}$ training trajectories is computed with parameters 
$$(\Delta t, n_{\text{train}}, L, L_{\text{train}}) = (.001, 100\times 100, 20, 20).$$

For testing, the reference solution is computed up to \( T = 1.6\) using a resolution of $(n_{\text{test}}, \Delta t) = (100\times 100, 0.001)$ with
\[(x_0, y_0, \alpha, \beta) = (1.032833,.034137,1.004777, .106782).\]
As in the 1D Burgers' equation case, only {\em smooth} solution profiles prior to the shock formation are available for training.

\subsubsection{2D KPP equation}\label{sub: 2D Burgers} 
As a final test case we consider the more challenging 2D KPP problem \cite{kpp}
\begin{equation}
    u_t + \big( \cos(u) \big)_x + \big( \sin(u) \big)_y = 0, \quad (x, y) \in \left[ -2, 2 \right]\times\left[ -2, 2 \right], \quad t > 0, 
    \label{eq: kpp 2d}
\end{equation}
subject to Dirichlet boundary conditions and initial conditions given by
\begin{equation}   
u\left( x, y, 0 \right) = \left\{
\begin{aligned}
&b e^{\frac{-r^2}{c^2}}, &\quad r>c\\
&b, & r < c
\end{aligned}.
\right.
    \label{eq: kpp 2d initial}
\end{equation}
where
\[\begin{+array}{ll}
c = 0.7 + \omega_c, & b = a + 1 + \cos 2 \pi x \sin 2 \pi y + \pi(1 + \cos 4 \pi x \sin 6 \pi y)\omega_b,\\
 a = 3.25\pi + 2\pi \omega_a,  &r = \sqrt{(x - x_0)^2 + (y - y_0)^2}. \\
\end{+array}\vspace{-0.2cm}\]

For training, the parameters $(x_0,y_0,\omega_a, \omega_b,\omega_c)$ in \eqref{eq: kpp 2d initial} are sampled independently from the uniform distribution \({\mathcal U}[-.25,.25]\). Each of the $N_{\text{traj}}$ training trajectories is computed using PyClaw with the HLLE Riemann solver, under the configuration $(\Delta t, n_{\text{train}}, L, L_{\text{train}}) = (.001, 100\times 100, 20, 20)$.

For testing, we fix the parameters in \eqref{eq: kpp 2d initial} to \(x_0= y_0=\omega_a=\omega_b=\omega_c = 0\). The reference solution is obtained on a spatial grid with $n_{\text{test}} = 100\times 100$, using a time step of $\Delta t = .001$ and integrated up to \( T = 0.6\).

\subsection{Observations with varying noise levels}
\label{subsec: noisy}
We now consider a training environment corrupted by noise, which is simulated  by applying additive noise to the  training trajectories corresponding to the PDE models introduced in \Cref{sec: pdeexamples}. 

\subsubsection{1D Burgers' equation} 
\label{subsubsec: 1D Burgers Noise}
Following the setup described in \Cref{sub:1DBurgers}, we generate perturbed training data 
\begin{equation}
\tilde{u}\big( x_j, t_l^{(k)} \big) = u\big( x_j, t_l^{(k)} \big) + \xi \overline{|u\left( x,t \right)|} \zeta_{j,l}, \quad k = 1,\dots,N_{\text{traj}}.
\label{eq: noise burgers}
\end{equation}
Here, \(\xi \in [0,1]\) is defined as the noise coefficient controlling the intensity of the perturbation while $\zeta_{j,l} \sim \mathfrak{N}\left( 0, 1 \right)$ denotes independent standard normal random variables, with indices $j = 1,\dots, n_{\text{train}}$ and $l = 0,\dots, L_{\text{train}}$.\footnote{Recall in this example we take $L_{\text{train}} = L$ (see \cref{eq:totalperiod} for the definition of $L$).} The quantity \(\overline{|u\left( x,t \right)|}\) represents the mean absolute value of exact solution \(u\left(x,t \right)\) over the entire dataset.  

\begin{figure}[h!]
        \includegraphics[width=\textwidth]{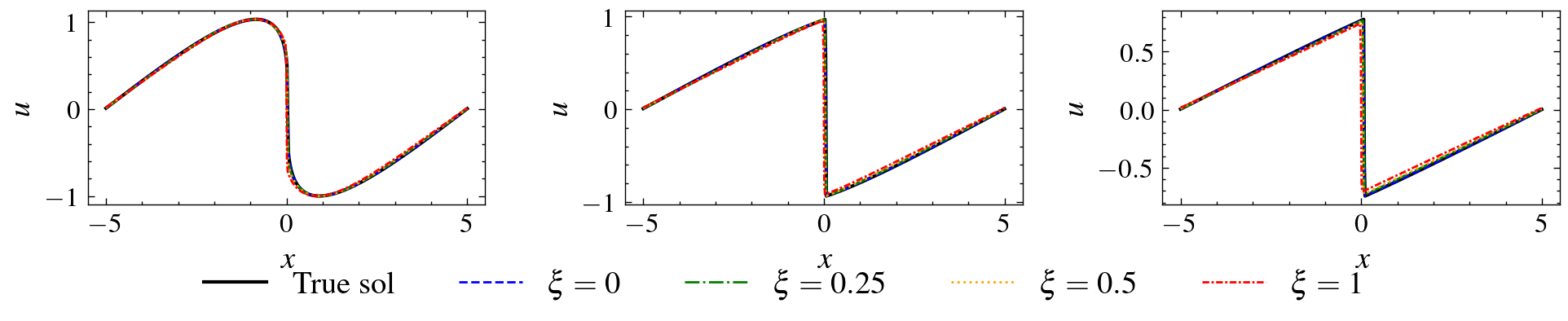}
    \caption{Comparison of the reference solution (black solid line) to the SymCLaw solution for 1D Burgers' equation at (left) \(t = 1\) (middle) \(t = 2\) (right) \(t = 3\) for $\xi = 0,.25,.5,1$ in  \eqref{eq: noise burgers}.}
    \label{fig: 1D Burgers Predictions Noise}
\end{figure}

\Cref{fig: 1D Burgers Predictions Noise} displays temporal snapshots of the SymCLaw prediction of the 1D Burgers' equation with noisy data, as defined in \eqref{eq: noise burgers}, for noise levels $\xi = 0,.25, .5$ and $1$. Remarkably, the model accurately captures the shock formation at \( t = 1.0\), {\em despite} being trained only on the interval $\mathcal{D}_{\text{train}} = [0, .1]$, during which the solution remains smooth. This extrapolatory capability was also observed in KT-ESCFN \cite{liu2024entropy} and NESCFN \cite{liu2025neural}. In comparison with these earlier methods, however, the SymCLaw solution demonstrates superior shock resolution in high-noise regimes. This improvement is anticipated, as shock speed is governed by the eigenvalues of the flux Jacobian, which are directly embedded into the SymCLaw.

\begin{figure}[h!]
        \centering
        \includegraphics[width=0.7\textwidth]{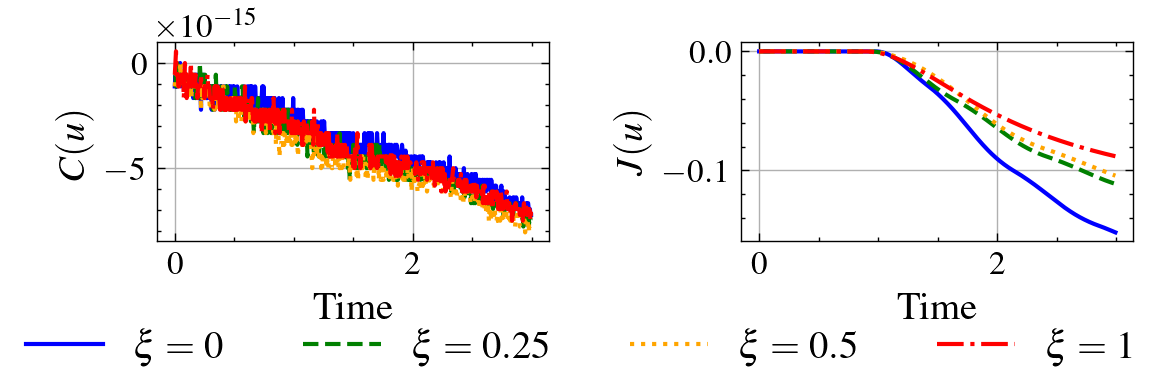}
       \caption{The discrete conserved quantity remainder $\mathcal{C}(u)$ defined in \eqref{eq:conservemetric} (left) and the discrete entropy remainder $\mathcal{J}(u)$ defined in \eqref{eq: discrete entropy} (right) for the SymCLaw prediction of the 1D Burgers' equation over the time interval $t \in [0,3]$ with noise levels $\xi = 0, .25, .5, 1$ in  \eqref{eq: noise burgers}.} 
    \label{fig: 1D Burgers Conservations Noise}
\end{figure}

\Cref{fig: 1D Burgers Conservations Noise} shows the time evolution of the discrete conservation remainder $\mathcal{C}({u})$, defined in \eqref{eq:conservemetric}, and the discrete entropy remainder $\mathcal{J}({u})$, defined in \eqref{eq: discrete entropy}, when training data are perturbed by different noise levels $\xi$. We observe that conservation in the SymCLaw prediction of the Burgers' equation is maintained at the level of $10^{-15}$, even at $T = 3$, well beyond the training interval $[0,.1]$. In addition, the discrete entropy remainders become clearly separated after shock formation, indicating that noise in the training data strongly influences the learned entropy function. Since training data are generated before shock formation, added noise at different spatial grids is interpreted as discontinuities, which in turn alters how the training algorithm handles the diffusion term. In principle, smooth profiles should introduce only minimal entropy dissipation, a point that will be further substantiated in the shallow water experiments. Overall these results demonstrate that SymCLaw yields consistent and robust long-term predictions for the 1D Burgers' equation, even when the training data are noisy. 

\subsubsection{Shallow water equation} 
\label{subsubsec: Shallow Water Noise}

We now examine the impact of noisy training data for the shallow water equation \cref{eq: shallow water equation}. As in the 1D Burgers’ case, zero-mean Gaussian noise is added to the training data within the domain $\mathcal{D}^{(k)}_{\text{train}}$, defined in \cref{eq:trainingperiod} for $k = 1,\dots,N_{\text{traj}}$. The perturbed data follows
\begin{equation}
    \tilde{\bf a}(x_j,t_l^{(k)})
     = {\bf a}(x_j,t_l^{(k)})
     + \xi\overline{|\bm{a}|}
       \zeta_{i,l},
       \label{eq:noise_shallow}
\end{equation}
where \({\bf a} = [h,hu]^\top\) denotes the vector of physical variables, \(\zeta_{j,l} \sim \mathfrak{N}({\bf 0}, \mathbb{I}_2)\) is a 2-dimensional standard normal vector, \( j = 1,\dots, n_{\text{train}}, l = 0,\dots, L_{\text{train}}\)\footnote{Recall that $L_{\text{train}} = L$ in this example.}, and \(\overline{|\bm{a}|}\) is the mean absolute value of the training data over the entire dataset.  We consider  noise intensity coefficient  \(\xi = 0, .25, .5,\) and $1$.

\begin{figure}[h!]
        \includegraphics[width=\textwidth]{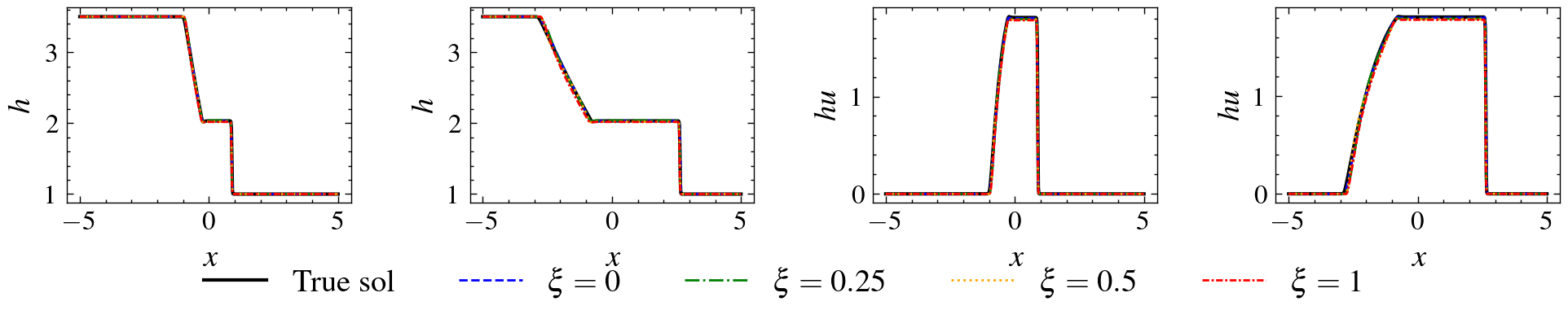}
        \caption{Comparison of the reference solution (black solid line) of height \(h\) and momentum \(hu\) for the SymCLaw of the shallow water equation with \(\xi = 0, .25, .5, 1\) in \eqref{eq:noise_shallow}: (left) \(t = .5\) of \(h\), (middle-left) \( t= 1.5\) of \(h\), (middle-right) \(t = .5\) of \(hu\), (right) \(t = 1.5\) of \(hu\).} 
    \label{fig: Shallow Water Predictions Noise}
\end{figure}

\Cref{fig: Shallow Water Predictions Noise} displays the predicted height \(h\) and momentum \(hu\) at \(t = .5,\) and \(t = 1.5\) across different noise levels $\xi$. The results show that the SymCLaw solution maintains accurate resolution of the shock structure even at the highest noise level $\xi = 1$. In contrast to the 1D Burgers' case, this robustness can be attributed to the inclusion of discontinuous states in the training data, which enables the model to directly learn shock-relevant features.

\begin{figure}[h!]
        \centering
        \includegraphics[width=0.9\textwidth]{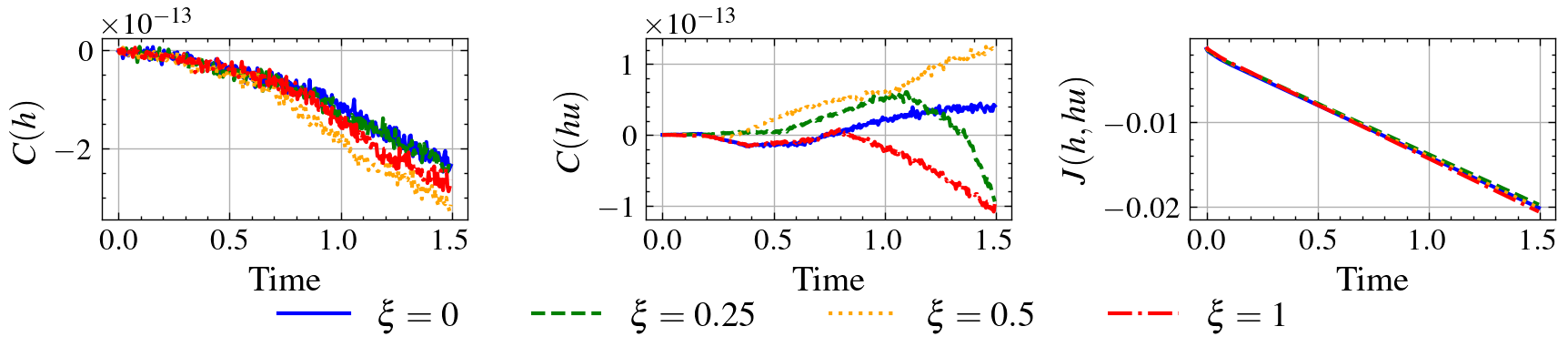}
    \caption{Discrete conserved quantity remainder  \eqref{eq:conservemetric} of height $\mathcal{C}(h)$ (left), the momentum \(\mathcal{C}(hu)\) (middle), and the discrete entropy remainder $\mathcal{J}([h,hu]^\top)$ (right) for the SymCLaw of the shallow water equations with noise levels $\xi = 0,.25,.5, 1$ in  \eqref{eq:noise_shallow}.} 
    \label{fig: Shallow Water Conservations Noise}
\end{figure}

\Cref{fig: Shallow Water Conservations Noise} (left-middle) presents the evolution of the discrete conservation remainders $\mathcal{C}(h)$ and \(\mathcal{C}(hu)\), while \Cref{fig: Shallow Water Conservations Noise} (right) depicts the corresponding discrete entropy remainder $\mathcal{J}([h,hu]^\top)$. The results demonstrate that conservation is preserved to within $\mathcal{O}(10^{-13})$, and that the entropy remainder remains non-positive, thereby confirming the entropy stability of the proposed SymCLaw method. Moreover, the consistent behavior of the entropy remainder in the shallow water case suggests that the entropy function is well trained, due to the presence of shock data in the training set.\footnote{Here, ``well trained" indicates that the proposed method achieves consistency in the shallow water case under noisy conditions; it does not imply that the correct entropy function for the shallow water equations can be recovered without additional regularization.}

\subsubsection{Euler's equations} 
\label{subsubsec: Euler Noise}
For the last 1D example, we investigate the impact of noise in the training data for  Euler's equations \cref{eq: Euler's equation}. As before, zero-mean Gaussian noise is added to the training data
\begin{equation}
    \tilde{\bf a}(x_j,t_l^{(k)})
     = {\bf a}(x_j,t_l^{(k)})
     + \xi\overline{|\bm{a}|}
       \zeta_{j,l}, \quad k = 1,2,\cdots, N_{\text{traj}};\quad  \xi \in \{0,.25,.5,1\},
       \label{eq:noise_euler}
\end{equation}
where \(\bm{a} = \left[ \rho, \rho u, E \right]^{T}\) represents the state vector, \(\zeta_{j,l} \sim \mathfrak{N}({\bf 0}, \mathbb{I}_3)\) is a standard 3-dimensional Gaussian vector, \(j = 1, \ldots, n_{\text{train}}\),  \(l = 0, \ldots,L_{\text{train}}\), with \(L_{\text{train}} < L\) (see Section \ref{sub:Eulers}). 
\begin{figure}[h!]
    \includegraphics[width=\textwidth]{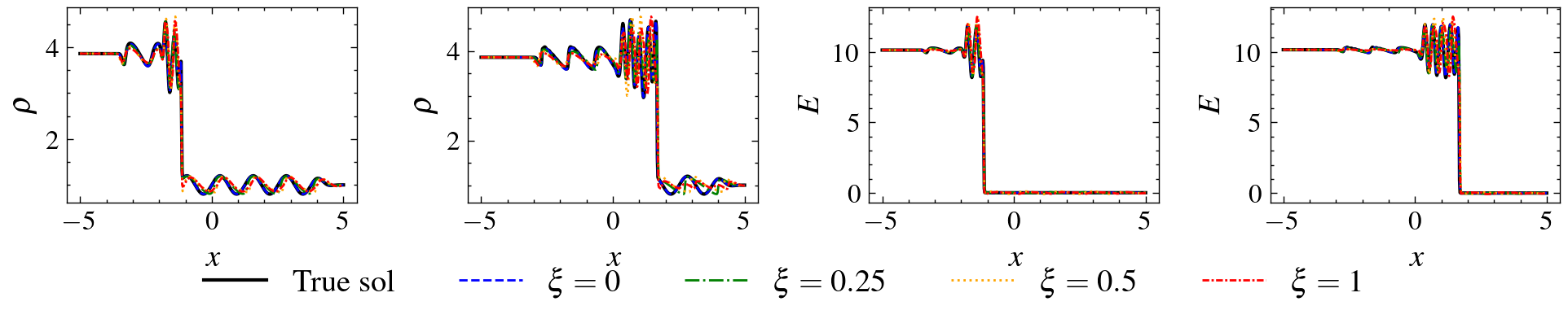}
        \caption{Comparison of the reference solution (black solid line) of the density \(\rho\) and energy \(E\) for the SymCLaw of the Euler's equation with \(\xi = 0, .25, .5, 1\) in  \eqref{eq:noise_euler}: (left) \(t = .8\) of \(\rho\), (middle-left) \(t = 1.6\) of \(\rho\), (middle-right) \(t =.8\) of \(E\), (right) \(t = 1.6\) of \(E\). }
    \label{fig: Euler Predictions Noise}
\end{figure}

\Cref{fig: Euler Predictions Noise} shows the model predictions for the density \(\rho\) and energy \(E\) across different noise levels $\xi$. We omit visualizations of the momentum \(\rho u\) as its behavior is qualitatively similar to that of \(E\). The prediction accuracy is comparable to that obtained by the KT-ESCFN (see \cite[section 5.2.3, Fig. 5.6]{liu2024entropy}) and the NESCFN (see \cite[section 5.2.3, Fig. 5.7]{liu2025neural}). It is worth emphasizing that noise in the training data strongly impacts the learning of both the flux potential \(\phi_{\bmu}(\bu)\) and the entropy function \(\eta_{\btheta}(\bu)\) in the entropy-stable flux \eqref{eq: Param entropy stable flux}, since their gradients and Hessians appear explicitly in the dissipative term. This sensitivity also explains the phase shift observed in the solution at later times when noisy data are used.

\begin{figure}[h!]
        \centering
        \includegraphics[width=0.9\textwidth]{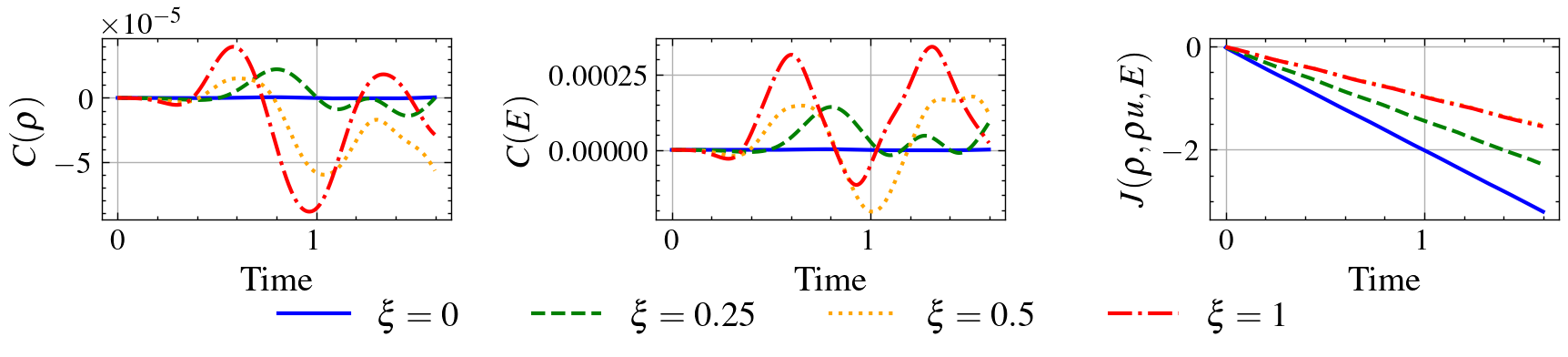}
   \caption{Discrete conserved quantity remainder \(\mathcal{C}(\rho)\)(left), \(\mathcal{C}(E)\)(middle), defined in \eqref{eq:conservemetric}, and discrete entropy remainder $\mathcal{J}([\rho, \rho u, E]^\top)$(right) for the SymCLaw of the Euler's equation with noise levels \(\xi = 0, .25, .5, 1\) in  \eqref{eq:noise_euler}.}
        \label{fig: Euler Conservations Noise}
\end{figure}

We also present the evolution of the discrete conservation remainders \(\mathcal{C}(\rho)\) and \(\mathcal{C}(E)\) in \Cref{fig: Euler Conservations Noise} (left-middle). These results confirm that conservation is preserved with the remainder on the order of $10^{-5}$. The larger variation, relative to machine precision, arises from the flux contributions induced by the inhomogeneous Dirichlet boundary conditions and from the temporal discretization error; it does not indicate any issue with the spatial approximation. In addition, the discrete entropy remainders $\mathcal{J}([\rho, \rho u, E]^\top)$ are non-positive across all noise levels, as illustrated in \Cref{fig: Euler Conservations Noise} (right). The effect of noise on the entropy function can be further understood through the evolution of the entropy remainder, with higher noise levels reducing entropy dissipation, and in turn causing a smaller contribution from the dissipative term in \eqref{eq: Param entropy stable flux}. 

\label{subsec: 2D Burgers noise}
\begin{figure}[h!]
    \centering
    \includegraphics[width=\linewidth]{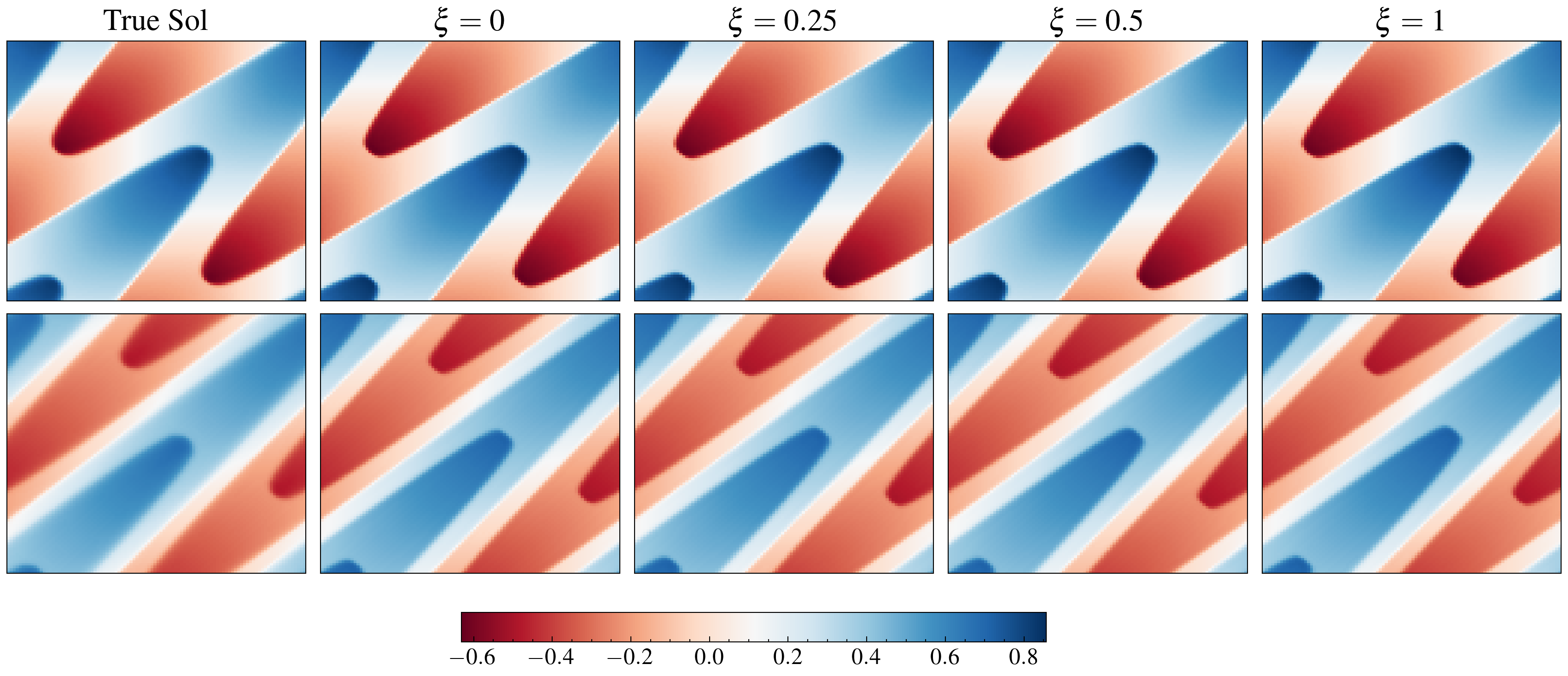}
    \caption{Comparison of the reference solution in 2D Burgers' equation \eqref{eq: Burgers' equation 2d} with the SymCLaw predictions for different noise levels $\xi = 0, .25,.5,1$. Top: \(t = 0.8\). Bottom: \(t = 1.6\).} 
    \label{fig: 2d Burgers Noise}
\end{figure}

\subsubsection{2D Burgers' equation}
We now consider the 2D Burgers' equation \eqref{eq: Burgers' equation 2d}. The noise is introduced following the same procedure as in \eqref{eq: noise burgers}, with noise levels \(\xi = 0, .25, .5, 1\). Recall that $L = L_{\text{train}} = 20$ with $\Delta t = .001$, and we again emphasize that the training trajectories only contain smooth profiles. \Cref{fig: 2d Burgers Noise} compares the reference solutions with the SymCLaw predictions  at $t = .8$ (top) and $1.6$ (bottom) for different noise levels. While the SymCLaw predictions are qualitatively accurate, even for $\xi = 1$, it is evident that it is difficult to capture the moving shock for \(\xi > .25\). 
This behavior is likely due to difficulties in accurately learning the local maximum wave speed, particularly under noisy and limited training data. Increasing the number of training trajectories and improving data quality may help the model to more precisely approximate the spectral radius. A more detailed investigation of these factors is left for future work. 

\begin{figure}[h!]
        \centering
        \includegraphics[width=0.7\textwidth]{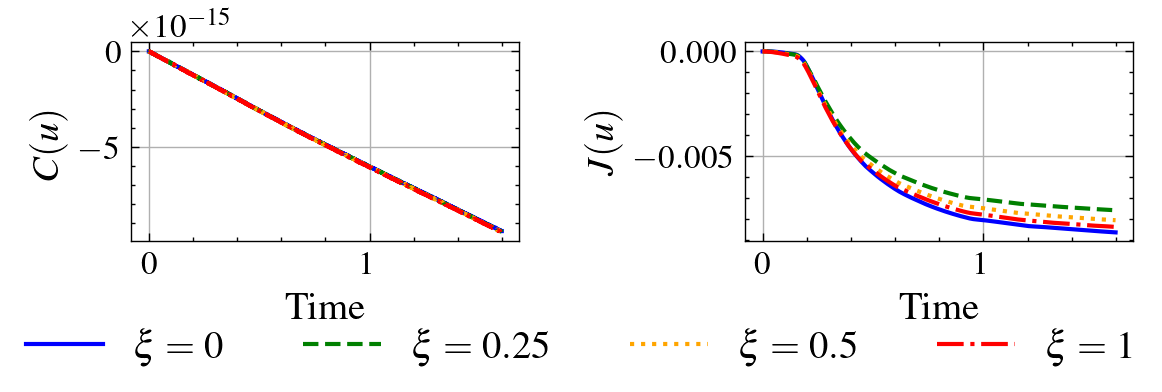}
    \caption{Discrete conserved quantity remainder  \eqref{eq:conservemetric} $\mathcal{C}(u)$ (left) and the discrete entropy remainder $\mathcal{J}(u)$ (right) for the parametric 2D Burgers' equations with noise levels $\xi = 0,.25,.5, 1$ in  \eqref{eq:noise_shallow}.} 
    \label{fig: 2D Burgers Conservations Noise}
\end{figure}

 \Cref{fig: 2D Burgers Conservations Noise} shows the evolution of the discrete conserved quantity remainder \(\mathcal{C}(u)\) and the discrete entropy remainder $\mathcal{J}(u)$, both of which behave similarly to the 1D case seen in \Cref{fig: 1D Burgers Conservations Noise}. 
The behavior of the discrete entropy remainder across different noise level indicates that SymCLaw solutions yield similar entropy diffusions, further demonstrating that the SymCLaw prediction maintains its robustness and entropy-stable properties in the multidimensional setting, even in considerable noisy data environments. 

\subsubsection{2D KPP equation} 
\label{subsec: 2D KPP noise}

\begin{figure}[h!]
    \centering
    \includegraphics[width=\linewidth]{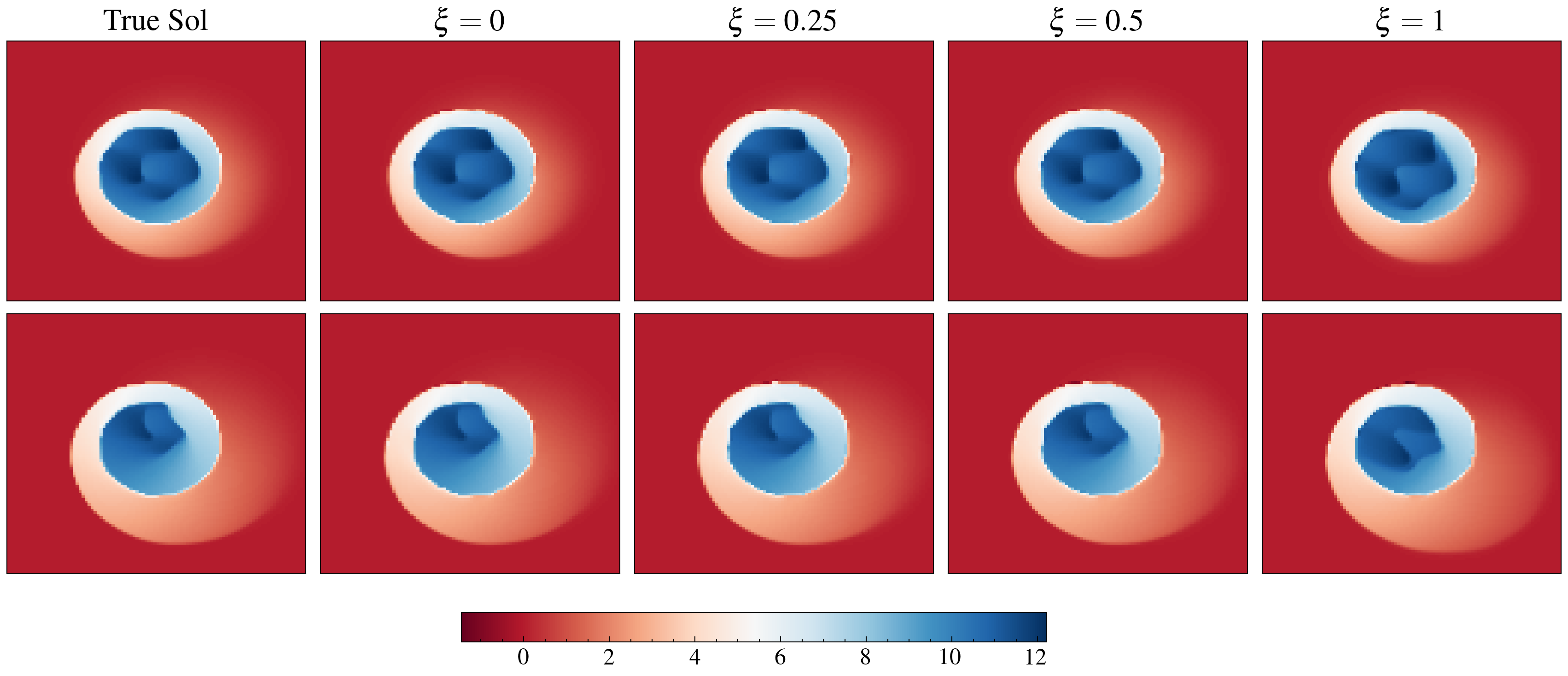}
    \caption{Comparison of the reference solution in 2D KPP \eqref{eq: kpp 2d} with the SymCLaw predictions for different noise levels $\xi = 0, .25,.5,1$. Top: \(t = 0.3\). Bottom: \(t = 0.6\).}
    \label{fig: 2d KPP Noise}
\end{figure}

We conclude by studying the recovery of the dynamics given noisy observations of the 2D KPP equation in \eqref{eq: kpp 2d}. Noise is introduced using the same procedure as in \eqref{eq: noise burgers}, with  \(\xi = 0, .25, .5, 1\). \Cref{fig: 2d KPP Noise} compares the reference solutions with the SymCLaw predictions for the 2D KPP equation at $t = .3$ (top) and $t = 0.6$ (bottom) across different noise levels. The results indicate that SymCLaw solutions retain good predictive accuracy for noise levels \(\xi = 0, .25, .5\).  We observe that the SymCLaw solution overestimates the shock speed when $u > 10$ in the case where $\xi = 1.0$. 
In contrast to the shallow water equations discussed in \Cref{subsubsec: Shallow Water Noise}, here the dynamics are successfully recovered for $\xi = 0;.25;.5$, even at the later time $t = 0.6$.  On the other hand, when $\xi = 1,$ the recovery is much less accurate.  This is despite the fact that shock information is available in the training data.
\begin{figure}[h!]
        \centering
        \includegraphics[width=0.7\textwidth]{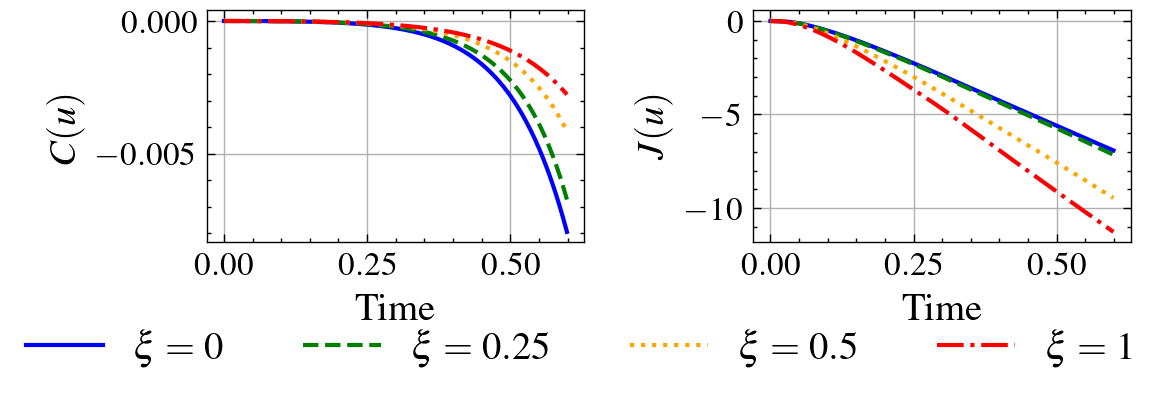}
    \caption{Discrete conserved quantity remainder $\mathcal{C}(u)$ (left) and the discrete entropy quantity remainder \(\mathcal{J}(u)\) (right) for the SymCLaw of the 2D KPP equation with noise levels $\xi = 0,.25,.5, 1$.} 
    \label{fig: 2D KPP Conservations Noise}
\end{figure}

Finally, \Cref{fig: 2D KPP Conservations Noise} presents the evolution of the discrete conserved quantity remainder \(\mathcal{C}(u)\) and the discrete entropy remainder $\mathcal{J}(u)$ for the SymCLaw solution of the 2D KPP equation. Unlike what is observed in  \Cref{fig: Euler Conservations Noise}, the discrete entropy remainder is much more sensitive to additional noise in the training data,  likely due to the reduced accuracy near the tails or the center. 

\section{Concluding remarks}\label{sec: conclusion}

This paper introduced the parametric hyperbolic conservation law (SymCLaw) method that simultaneously enforces conservation, entropy stability, and hyperbolicity within a unified learning paradigm. We further proposed an entropy-stable numerical flux scheme, enabling the framework to be flexibly integrated with existing classical numerical methods. Numerical experiments demonstrate that SymCLaw solutions maintain stability and accuracy even under noisy training conditions, generalize effectively to unseen initial conditions, and preserve the entropy inequality throughout long-term predictions.

Future work may proceed along several important directions. A natural extension is the incorporation of high-order discretizations, such as discontinuous Galerkin methods, to address both higher-dimensional problems and complex geometries. Establishing an accurate mechanism for  inferring boundary conditions is also an important direction, as is the development of rigorous theoretical guarantees, including convergence rates, error bounds, and long-time stability analysis, to complement the empirical findings. Beyond the deterministic setting, adapting the framework to uncertainty quantification and stochastic conservation laws would broaden SymCLaw's applicability to real-world systems with inherent randomness. 

\section*{Acknowledgements} 
This work was partially supported by the DOD (ONR MURI) grant \#N00014-20-1-2595 (Liu and Gelb), the DOE ASCR grant \#DE-SC0025555 (Gelb), the NSF grants \#DMS-2231482 (Zhang) and \#DMS-2513924 (Zhang), Simons Foundation through travel support for mathematicians (Zhang), and Ken Kennedy Institute at Rice University (Zhang).

\appendix 

\small
\bibliographystyle{siamplain}
\bibliography{reference}

\end{document}